%% file: MSOM_submission.tex
\documentclass[msom,nonblindrev]{informs3}
\usepackage{amsmath,amssymb,amsfonts}
\usepackage{bm}
\usepackage{booktabs}
\usepackage{graphicx}
\usepackage{xcolor}
\usepackage[normalem]{ulem}
\usepackage{hyperref}
\usepackage{xurl}
\usepackage{natbib}
\usepackage{microtype}
\usepackage{multirow}
\hypersetup{
  colorlinks=true,
  linkcolor=blue,
  filecolor=blue,
  urlcolor=blue,
  citecolor=blue,
}
\DoubleSpacedXI
\TheoremsNumberedThrough
\ECRepeatTheorems
\EquationsNumberedThrough

\bibpunct[]{(}{)}{,}{a}{}{,}
\def\bibfont{\small}

\OneAndAHalfSpacedXI
\begin{document}

% \RUNAUTHOR{Anonymous Authors}

\RUNTITLE{DRO Dynamic Inventory Control: Automatic Recovery under Martingale Demand}
\TITLE{Distributionally Robust Dynamic Inventory Control: Automatic Recovery under Martingale Demand}

\ARTICLEAUTHORS{
\AUTHOR{Jianzhi Leng$^{1,ab}$, Daniel Zhuoyu Long$^{2,c}$, Lun Ran$^{3,a}$, Chi Xie$^{4,a}$}
\AFF{$^{a}$School of Management, Beijing Institute of Technology, Beijing, China.\\
$^{b}$Digital Economy and Policy Intelligentization Key Laboratory of Ministry of Industry and Information Technology.\\
$^{c}$Department of Systems Engineering and Engineering Management, The Chinese University of Hong Kong, Hong Kong, China.
\\\EMAIL{$^{1}$lengjianzhi@bit.edu.cn, 
$^{2}$zylong@cuhk.edu.hk,
$^{3}$ranlun@bit.edu.cn, $^{4}$xc@bit.edu.cn}}
}

\RUNAUTHOR{Leng et al.}
\ABSTRACT{
\textbf{Problem definition:}
We study dynamic inventory control for battery-swapping
stations that combine central replenishment with automatic
local recovery. Central replenishment improves immediate
battery availability but reduces the depleted-battery pool
available for subsequent recovery, creating an intertemporal
trade-off. Managing this trade-off is particularly challenging
when historical demand data are limited and future demand
is uncertain and temporally dependent.
\textbf{Methodology/results:}
We develop a time-consistent distributionally robust dynamic
program under martingale demand ambiguity, using the latest
observed demand as the next conditional mean without requiring
a fully specified demand distribution. We establish an optimal
order-up-to policy for every finite horizon and identify two
critical recovery rates governing replenishment decisions.
Under a regular penalty condition, we derive explicit
finite-horizon policies and worst-case demand distributions,
showing that sufficiently strong automatic recovery makes
central replenishment unnecessary in nonterminal periods.
Extensions to submartingale and supermartingale demand
reveal contrasting behaviors, including eventual
full-capacity replenishment under upward ambiguity and
recovery-induced positive replenishment under downward
ambiguity.
\textbf{Managerial implications:}
Numerical experiments demonstrate that the robust policy
is particularly competitive when demand is highly dispersed
and historical information is limited, whereas its advantage
diminishes when demand becomes more predictable or the
conditional-mean assumption becomes less reliable.
Our findings also reveal that greater reliance on local
recovery may compromise service availability, highlighting
the importance of jointly considering demand information,
recovery capability, and service requirements when making
replenishment decisions.
}
\KEYWORDS{distributionally robust optimization; dynamic inventory control;
martingale demand; battery swapping; automatic recovery}

\maketitle

\section{Introduction}
\begin{figure}[ht]
    \centering
    \includegraphics[width=0.68\linewidth]{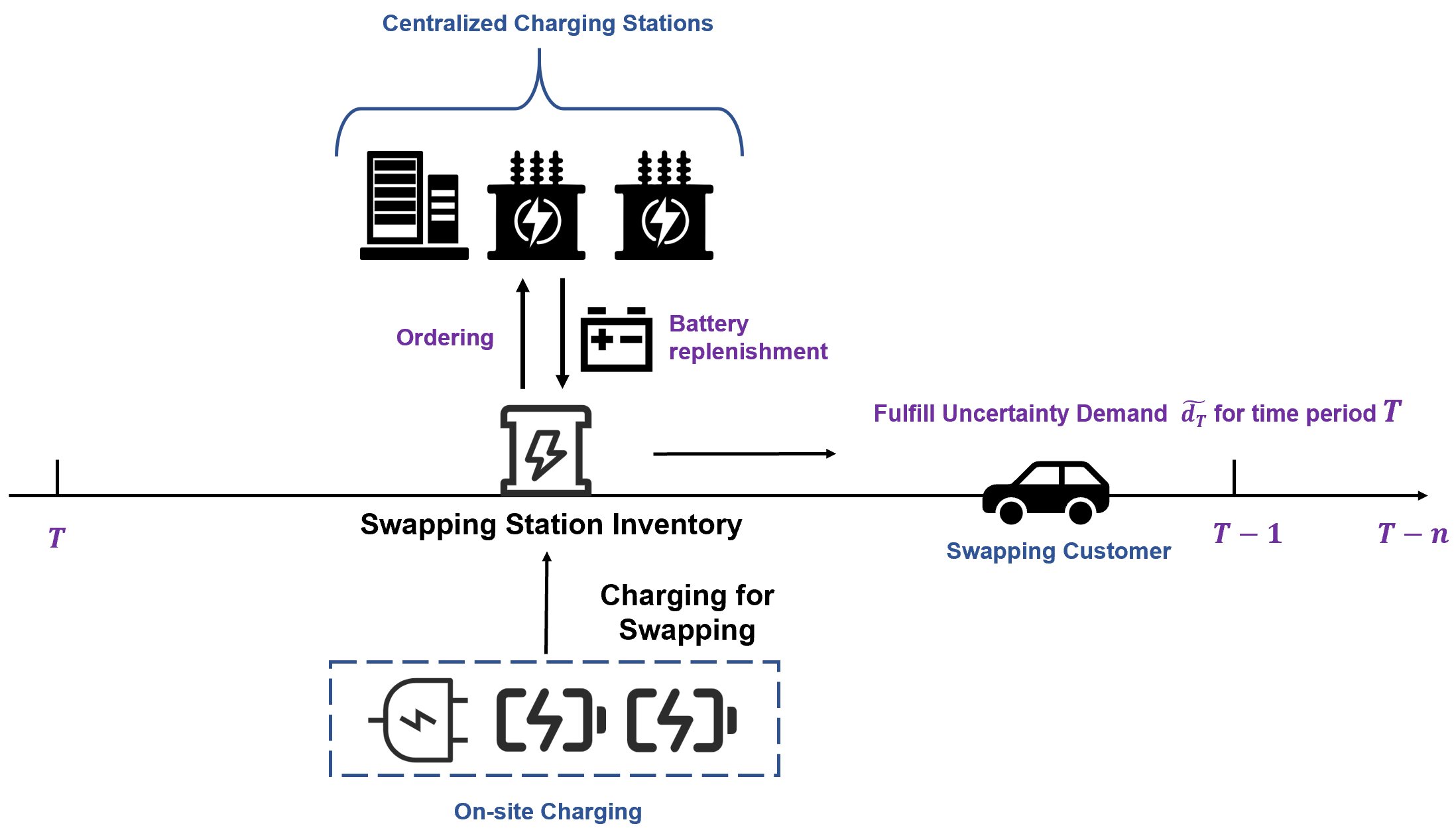}
    \caption{Schematic of the SLCA operating system}
    \label{fig:SLCA}
\end{figure}

The rapid growth of electric vehicles (EVs) has increased the need for
reliable charging and battery-swapping services. The International Energy Agency projects that the global EV fleet could reach about 250 million vehicles by 2030 under its Stated Policies Scenario, excluding two- and three-wheelers (\citealt{IEA_GlobalEVOutlook2025}). Battery swapping offers faster energy replenishment than plug-in charging, but its infrastructure deployment has lagged behind EV adoption (\citealt{patel2024impact}, \citealt{chen2024charging}). A central operational challenge is to maintain sufficient ready-battery inventory despite limited station capacity, volatile demand, and recharging delays. One way to improve battery availability is to combine centralized replenishment with local charging. We refer to this hybrid configuration as ``Swap-Locally, Charge-Adaptively (SLCA)'' as illustrated in Figure~\ref{fig:SLCA}. Ready batteries can be supplied from a centralized charging facility, while depleted batteries retained at the swapping station can be recharged locally and returned to the ready-battery inventory. Related centralized charging-and-swapping configurations have been studied in the literature (\citealt{liu2018distributed}, \citealt{li2022battery},
\citealt{qi2023scaling}). Charging and battery-swapping infrastructure has also received policy support in cities such as Shanghai and Shenzhen (\citealt{Shanghai2022}, \citealt{Shenzhen2023}).

Our operational model captures a distinctive feature of this hybrid
system: central replenishment is actively controlled, whereas local
charging provides automatic recovery. After the station chooses its
ready-battery level, a fixed fraction of the depleted batteries remaining on site is recharged before the next decision epoch. Replenishing more ready batteries therefore reduces immediate shortage exposure but leaves fewer depleted batteries available for local recovery, creating an intertemporal trade-off between current service and future endogenous supply. This trade-off becomes especially difficult to manage when demand is serially dependent but its probability law is not reliably known. A central question is therefore how the recovery intensity changes the value of centralized replenishment as the remaining horizon evolves. Although correlated demand and evolving forecasts have long been studied in dynamic inventory models (\citealt{scarf1959bayes}, \citealt{scarf1960some}, \citealt{chen2000quantifying}, \citealt{gallego2001integrating},
\citealt{dong2003optimal}, \citealt{dehaybe2024deep}), these approaches
typically require a specified stochastic demand model. Such a model may
be difficult to estimate and sensitive to misspecification (\citealt{badinelli1990inventory}, \citealt{fildes1995impact},
\citealt{fildes2008forecasting}), particularly for newly deployed
swapping stations with limited historical data.

Distributionally robust optimization (DRO) provides a natural way to
address this form of model uncertainty by optimizing against the
worst-case distribution within an ambiguity set constructed from
available information (\citealt{scarf1957min},
\citealt{bertsimas2010models}, \citealt{delage2010distributionally}).
Such sets may be defined through moments
(\citealt{scarf1957min}, \citealt{hanasusanto2015distributionally},
\citealt{raad2025hybrid}) or through distances from a reference
distribution (\citealt{gao2024wasserstein}, \citealt{kim2024data}).
In multistage problems, however, a joint ambiguity set alone does not
generally provide a dynamically consistent model. Distributional
restrictions may couple uncertainty across histories, preventing the
problem from decomposing after new demand is observed and causing
continuation decisions obtained by reoptimization to differ from those
implied by the original plan, i.e., time inconsistency
(\citealt{sayman2009investigation}, \citealt{xin2021time}). This
motivates ambiguity models that specify how admissible conditional
distributions evolve with the observed demand history.

Several studies develop dynamic programming approaches for
distributionally robust multistage problems
(\citealt{adida2006robust}, \citealt{pichler2021mathematical},
\citealt{maggioni2025sampling}). A key requirement for such formulations
is \textit{rectangularity}: after new information is observed, the
admissible joint distributions can be decomposed into stagewise
conditional ambiguity sets
(\citealt{iyengar2005robust}, \citealt{shapiro2016rectangular},
\citealt{xin2022distributionally}). This structure ensures time
consistency and allows the worst-case expectation to be evaluated
recursively through dynamic programming. Even under rectangularity,
however, multistage inventory problems remain difficult to analyze, and
most existing approaches rely on numerical methods. Although structural
policies such as nonstationary or history-dependent $(s,S)$ rules are
known for broad classes of robust inventory problems
(\citealt{ahmed2007coherent}, \citealt{klabjan2013robust}), explicit
finite-horizon policies under dynamically dependent demand remain scarce.

We therefore model demand through stagewise conditional-mean information.
Specifically, the ambiguity set is constructed so that, after observing
the latest demand, the next conditional distribution may be any
distribution on the bounded support whose mean equals that observation.
Equivalently, the demand process satisfies the martingale condition
\begin{equation*}
\mathbb E[\widetilde D_t\mid\mathcal F_{t+1}]
=
\widetilde D_{t+1},
\qquad t\in[T],
\end{equation*}
with $\widetilde D_{T+1}:=\mu$. Because the ambiguity set is specified
through these history-by-history conditional moment restrictions, it
admits the rectangular decomposition required for dynamic programming. The martingale restriction leaves the conditional law otherwise unspecified beyond its conditional mean, thereby capturing serial dependence without requiring a fully specified transition law. This makes it a parsimonious short-horizon demand model for settings in which historical data are limited but demand information is updated sequentially.

We make three main contributions.
\begin{itemize}

\item We develop a time-consistent distributionally robust
inventory model that integrates central replenishment
with automatic local recovery under martingale demand
ambiguity. Using only sequentially updated conditional-mean
information, the model accommodates temporally dependent
demand without requiring a fully specified distribution.
We establish the existence of an optimal order-up-to
policy for every finite horizon.

\item We characterize how automatic recovery shapes optimal
replenishment decisions. Two critical recovery rates
identify four operating regimes and determine how
replenishment incentives evolve with the planning horizon.
We derive exact two-period solutions and, under a regular
penalty condition, explicit finite-horizon policies and
worst-case demand distributions. In particular, central
replenishment can remain beneficial at low recovery rates
but becomes unnecessary in all nonterminal periods beyond
a critical threshold. We also provide approximation
policies with bounded optimality gaps for the complementary
case.

\item We investigate how directional demand dependence
changes replenishment behavior. Under submartingale
ambiguity, sufficiently long horizons can make full-capacity
replenishment optimal, whereas supermartingale ambiguity
reveals that stronger local recovery can sometimes
encourage central replenishment. Numerical experiments
further demonstrate the value and limitations of M-DRO:
it is particularly competitive under highly dispersed
demand and limited historical information, while its
advantage depends on the reliability of conditional-mean
information and the operating environment.

\end{itemize}

Section~\ref{literature} reviews related literature,
and Section~\ref{Model} presents the model.
Sections~\ref{sec:two-period}--\ref{sec:directional-drift}
develop the analytical results, followed by numerical
experiments in Section~\ref{Numerical} and conclusions
    in Section~\ref{sec:conclusion}. Proofs and supplementary
results appear in the Online Appendix.

\textbf{Notations:} For any positive integer $n$, let
$[n]:=\{1,\ldots,n\}$. Boldface is used for vectors, and tildes distinguish random demand from its realizations.
Let $\mathcal P_0(\mathcal Z)$ denote the set of
Borel probability measures supported on $\mathcal Z$.
Expectations under $\mathbb P$ and $\mathbb Q$ are
denoted by $\mathbb E_{\mathbb P}[\cdot]$ and
$\mathbb E_{\mathbb Q}[\cdot]$, respectively.
We use $\delta_z$ for the Dirac measure at $z$
and $(z)_+:=\max\{z,0\}$.

\section{Related Literature}
\label{literature}
Highly relevant literature can be divided into two streams: the operations management research of EV battery swapping stations, and the discussion of distributionally robust optimization in inventory management, both in static and dynamic.
\subsection{Battery-Swapping Inventory and Recovery Operations}
EVs with battery-swapping services have been experiencing a resurgence in recent years (\citealt{mak2013infrastructure}, \citealt{wu2021survey}, \citealt{zhang2025battery}). Battery-swapping research studies network design, inventory, and charging operations. For example, \cite{qi2023scaling} develop a joint location and repairable-inventory model and show that centralized charging can achieve significant cost efficiencies while maintaining service performance. \cite{liang2022real} study the real-time operation of battery swapping–charging systems and propose a multi-agent deep reinforcement learning framework. Their approach coordinates swapping and charging decisions under stochastic demand and system dynamics. \cite{shaker2023joint} investigate the joint planning of electric-vehicle battery swapping stations and the distribution grid under centralized charging. They develop an integrated planning model that captures the interactions between station deployment, battery charging, and grid constraints. These studies address different planning and operating decisions, including dynamic charging control. The single-station benchmark here abstracts from network allocation and charging schedules to isolate how automatic recovery changes replenishment geometry. To account for demand uncertainty, these works often formulate two-stage stochastic optimization models for scheduling and routing decisions (see \citealt{schneider2018optimization}, \citealt{sun2019robust}, \citealt{ni2020inventory}, \citealt{he2022sizing}). Some studies model inventory control problems as Markov decision processes (MDP). However, solving the resulting DP via backward induction often becomes computationally intractable due to the curse of dimensionality. Consequently, most of this literature resorts to approximate solution methods (see \citealt{asadi2021stochastic}, \citealt{asadi2022monotone}). Our work provides an explicit finite-horizon characterization of the optimal policy under the regular penalty condition. Its unit of analysis is a single station with immediately available central supply and a fixed recovery rule. Central-facility capacity, transportation lead times, multi-station allocation, and charging queues are outside this benchmark.

\subsection{Distributionally Robust Optimization in
Inventory Management}
A conventional approach to decision making under uncertainty is stochastic optimization. Common stochastic optimization methods include sample average approximation (\citealt{shapiro2002minimax}) and stochastic approximation methods (\citealt{shapiro2000rate}). Moreover, the distribution of random variables can be difficult to predict or fit accurately, and imposing a probability distribution can result in biased decisions or overfitting. DRO addresses this uncertainty by optimizing over a specified set of probability laws; its guarantees depend on the information used to construct that set. In practice, \cite{qiu2021distributionally} develop a distributionally robust optimization framework for multi-product inventory decisions under budget constraints, accounting for both demand and yield uncertainties. Their approach provides inventory policies that are robust to distributional ambiguity while controlling total investment cost. \cite{kang2015distribution} use distributional information to reduce conservatism in robust linear optimization while retaining probabilistic feasibility guarantees, and apply the method to inventory control with quality-of-service constraints. \cite{yin2024distributionally} study a distributionally robust multi-period humanitarian relief network design problem that integrates facility location, inventory, allocation, and evacuation decisions. Their model enhances network resilience by jointly accounting for uncertainty across multiple planning stages. For more details, \cite{zhang2025robusta} and \cite{zhang2025robustb} provide a comprehensive survey of robust optimization approaches in inventory management, reviewing modeling frameworks, uncertainty sets, and solution methods. These papers synthesize theoretical developments and applications, highlighting how robust optimization addresses demand, cost, and supply uncertainties in inventory systems.

For multistage decisions, the distinction between a joint ambiguity set and a rectangular family of conditional laws is central. Robust dynamic programming provides the decomposition for rectangular transition uncertainty (\citealt{iyengar2005robust}), while distributional constraints linking histories require separate attention to dynamic consistency (\citealt{shapiro2016rectangular}). We use this conditional framework to analyze the inventory coupling created by automatic local recovery. Our analysis proceeds from one-dimensional conditional moment problems
to demand envelopes, recursive demand points and inventory thresholds, and finite-horizon
closed-form policies.

\section{Model Formulation}
\label{Model}
We formulate a finite-horizon inventory problem with central replenishment, automatic local recovery, and rectangular conditional-mean ambiguity. The Bellman formulation is valid at every finite horizon; the explicit closed-form representations developed later require additional regime-specific conditions.

\subsection{System dynamics and demand information}

We consider $T$ operating periods. The index $t$ is the number of periods
remaining, so decisions occur in the order $T,T-1,\ldots,1$. At the beginning
of period $t$, the station has $y_t$ ready batteries and chooses a
post-replenishment level $x_t\in[y_t,U]$. In this stylized model, $U>0$
is both the ready-inventory capacity and the size of a fixed on-site battery
pool: when $x_t$ batteries are ready, $U-x_t$ battery equivalents are depleted
and eligible for local recovery. Central replenishment is immediately
available, equals $x_t-y_t$, and replaces the same number of depleted battery
equivalents in the on-site pool. The demand vector is
$\widetilde{\bm D}:=(\widetilde D_T,\ldots,\widetilde D_1)$, with
realization $\bm D:=(D_T,\ldots,D_1)$.

Unmet demand is lost and incurs a unit penalty $b$, excess ready inventory incurs a unit holding cost $h$, and central replenishment costs $c$ per battery (\citealt{dana2001note}, \citealt{shapiro2012minimax},
\citealt{wang2012multiordering}). Local charging provides automatic recovery. After $x_t$ is selected, $U-x_t$ depleted batteries remain in the charging pool; a fixed fraction $\alpha\in[0,1)$ completes charging before the next decision epoch at unit cost $w$. We exclude the degenerate endpoint \(\alpha=1\), at which the entire pre-demand depleted pool is recovered within one review period. Write $\beta:=1-\alpha$. The one-period cost function is
\begin{equation*}
C(y,x,D)
:=
c(x-y)+w\alpha(U-x)
+b(D-x)_+
+h(x-D)_+,
\end{equation*}
with the positive-part notation defined above.

\begin{figure}[ht]
    \centering
    \includegraphics[width=0.85\linewidth]{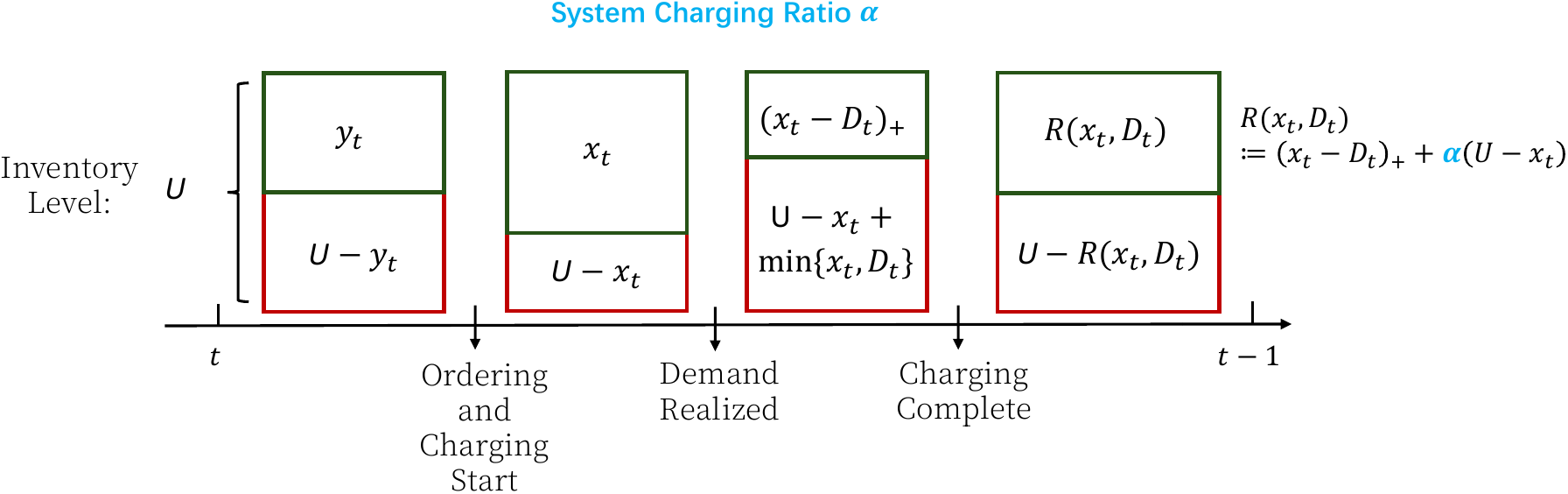}
    \caption{Sequence of Events}
    \label{fig:sequence-events}
\end{figure}

Figure~\ref{fig:sequence-events} summarizes the within-period sequence. Once
$x_t$ is chosen, the current recovery volume $\alpha(U-x_t)$ is fixed. Demand
$D_t$ is then realized and served from the $x_t$ ready batteries. Batteries
returned by swaps during period $t$ join later charging batches and do not
complete charging within the same review interval. The batch committed before
demand completes before the next decision epoch, giving
$
y_{t-1}
=
R(x_t,D_t)
:=
(x_t-D_t)_+ + \alpha(U-x_t)$. For every $(x,D)\in[0,U]^2$, $R(x,D)\in[0,U]$.
The second state variable records the conditional-demand information. At the
initial epoch, $\mu_T=\mu$ is given. For $t<T$, the latest observed demand is
$D_{t+1}$, so define $D_{T+1}:=\mu$ and $\mu_t:=D_{t+1}$. After observing
$D_t$, the next state is $(y_{t-1},\mu_{t-1})=\bigl(R(x_t,D_t),D_t\bigr)$. Thus $(y_t,\mu_t)$ contains the ready inventory and the reference value that
determines the next conditional ambiguity set.

Set $\widetilde D_{T+1}\equiv\mu$ and
$\mathcal F_{T+1}:=\{\emptyset,\Omega\}$. For $t\in[T]$, let $\mathcal F_t:=\sigma(\widetilde D_T,\ldots,\widetilde D_t)$. A policy $\pi=(\pi_T,\ldots,\pi_1)$ consists of deterministic measurable
rules based on the observed history, with
$\pi_t:[0,U]^{T-t}\to[0,U]$ and
$x_t^\pi=\pi_t(\widetilde D_T,\ldots,\widetilde D_{t+1})$. Feasibility requires $y_t^\pi\leq x_t^\pi\leq U$, and inventory evolves by
$y_{t-1}^\pi=R(x_t^\pi,\widetilde D_t)$.
For each initial inventory $y\in[0,U]$, let $\mathcal A(y)$ denote the set
of admissible policies satisfying these conditions with $y_T^\pi=y$
for every demand history.

We define the joint martingale ambiguity set as
\begin{equation*}
\mathcal{MAR}(\mu)
:=
\left\{
\mathbb P\in\mathcal P_0([0,U]^T):
\mathbb E_{\mathbb P}
[\widetilde D_t\mid\mathcal F_{t+1}]
= \widetilde D_{t+1},
\quad t\in[T]
\right\}.
\end{equation*}
For a one-period reference mean $m\in[0,U]$, define
\begin{equation*}
\mathcal M(m)
:=
\left\{
\mathbb Q\in\mathcal P_0([0,U]):
\mathbb E_{\mathbb Q}[\widetilde D]=m
\right\}.
\end{equation*}
Under any $\mathbb P\in\mathcal{MAR}(\mu)$, the conditional law of
$\widetilde D_t$ after history $(D_T,\ldots,D_{t+1})$ belongs to
$\mathcal M(D_{t+1})$; initially the law of $\widetilde D_T$ belongs to
$\mathcal M(\mu)$. Conversely, a sequence of measurable conditional kernels
chosen from these one-period sets composes to an admissible joint law. This stagewise representation provides the rectangular structure underlying the Bellman recursion.

For initial inventory $y$ and initial reference mean $\mu$, the robust control problem for $T$ periods is
\begin{equation*}
\mathcal J_T(y,\mu)
:=
\inf_{\pi\in\mathcal A(y)}\sup_{\mathbb P\in\mathcal{MAR}(\mu)}
\mathbb E_{\mathbb P}
\left[
\sum_{t=1}^{T}
C\!\left(y_t^\pi,x_t^\pi,\widetilde D_t\right)
\right],
\end{equation*}
where $y_T^\pi=y$ and
$y_{t-1}^\pi=R(x_t^\pi,\widetilde D_t)$.

\subsection{Bellman recursion and analytical foundations}
For each $t\in\{0,1,\ldots,T\}$, let
$\hat V^t(y,\mu)$ denote the optimal worst-case cost-to-go
when $t$ operating periods remain, the current ready inventory is $y$,
and the demand mean is $\mu$. The terminal condition
is $\hat V^0(y,\mu)=0$ on $[0,U]^2$, that is, there is no salvage value or post-horizon operating cost. Since the martingale ambiguity set is rectangular, the worst-case expected
cost can be evaluated recursively over admissible conditional distributions.
Hence, for $t\in[T]$, the cost-to-go function satisfies the Bellman recursion
\begin{equation*}
\begin{aligned}
\hat V^t(y,\mu)
=\min_{x\in[y,U]}\Bigg\{&
c(x-y)+w\alpha(U-x)\\
&+\sup_{\mathbb Q\in\mathcal{M}(\mu)}
\mathbb{E}_{\mathbb Q}
\left[
b(\widetilde D-x)_+
+h(x-\widetilde D)_+
+\hat V^{t-1}(R(x,\widetilde D),\widetilde D)
\right]
\Bigg\},\quad\forall t\in [T].
 \end{aligned}
\end{equation*}
After demand $D$, the continuation state is $(R(x,D),D)$; unmet demand is lost and does not enter the state. The following proposition formalizes the equivalence with the multistage minimax problem.
\begin{proposition}
\label{prop:dp-equivalence}
For every $(y,\mu)\in[0,U]^2$,
\begin{equation*}
    \mathcal J_T(y,\mu)=\hat V^T(y,\mu).
\end{equation*}
\end{proposition}
Proposition~\ref{prop:dp-equivalence} also establishes
attainment of the Bellman minimum and the conditional suprema. By rectangularity, the attaining conditional kernels can be concatenated into an attaining joint distribution in $\mathcal{MAR}(\mu)$. Together with the
Bellman equivalence, these attainment properties yield time consistency.
\begin{corollary}[Time consistency]
\label{cor:time-consistency}
There exists an optimal deterministic Markov policy
whose decision at each state $(y,\mu)$ is selected
from the Bellman minimizers. For every
$t\in[T]$ and $(y,\mu)\in[0,U]^2$, this policy
is optimal for the corresponding $t$-period
problem and attains $\hat V^t(y,\mu)$.
\end{corollary}

Under the conditional mean and support restrictions, the worst-case
expectation in the Bellman recursion can be expressed through a
concave envelope (\citealt{shapiro2012minimax}). Moreover, a
worst-case distribution can be chosen with support on at most
two demand points.
\begin{lemma}
\label{lem:moment-envelope}
For any continuous $\varphi:[0,U]\to\mathbb R$ and
$\mu\in[0,U]$, a worst-case distribution $\mathbb Q^*\in\arg\max_{\mathbb Q\in\mathcal M(\mu)}\mathbb E_{\mathbb Q}[\varphi(\widetilde D)]$
can be chosen with support on at most two contact points
$D^-\leq\mu\leq D^+$. For $D^-<D^+$, the corresponding probabilities are
\[
\mathbb Q^*(D^-)
=
\frac{D^+-\mu}{D^+-D^-},
\qquad
\mathbb Q^*(D^+)
=
\frac{\mu-D^-}{D^+-D^-}.
\]
If $D^-=D^+=\mu$, then $\mathbb Q^*=\delta_\mu$.
\end{lemma}

For a fixed post-replenishment level $x$, define the demand-dependent
cost, comprising current shortage and holding costs and the optimal
continuation cost, as
\begin{equation*}
F_t^x(D)
:=
b(D-x)_+ + h(x-D)_+
+\hat V^{t-1}\!\left(R(x,D),D\right),
\quad t\in[T].
\end{equation*}
By Lemma~\ref{lem:moment-envelope}, the worst-case expectation of
$F_t^x$ is attained by a distribution supported on at most two demand
points. Including the replenishment and recovery costs and separating
the state-dependent term $-cy$, we define
\begin{equation*}
\begin{aligned}
B_t(x,\mu)
:=\;&cx+w\alpha(U-x)\\
&+\sup_{\mathbb Q\in\mathcal M(\mu)}
\mathbb E_{\mathbb Q}[F_t^x(\widetilde D)],\quad t\in[T].
\end{aligned}
\end{equation*}
The Bellman equation can therefore be written as 
\begin{equation*}
\hat V^t(y,\mu)=-cy+\min_{x\in[y,U]}B_t(x,\mu),\quad t\in [T]. 
\end{equation*}
To characterize the replenishment policy analytically, we impose the
following strict ordering on the four unit costs. No ordering is imposed on the
recovery thresholds introduced below.

\begin{assumption}
\label{ass:cost-region}
$0<h<w<c<b$.
\end{assumption}

The ordering reflects inexpensive inventory carryover, local charging that is
cheaper than central replenishment, and a shortage penalty exceeding the
central replenishment cost. Define the two critical recovery rates
\begin{equation*}
\alpha_h:=\frac{h}{w},
\qquad
\alpha_b:=\frac{b-c}{c-w}.
\end{equation*}

As we show in the subsequent analysis, these two thresholds
play distinct roles: $\alpha_h$ governs how the recursive
coefficients evolve with the remaining horizon, whereas
$\alpha_b$ marks the boundary for the optimality of no
replenishment in nonterminal periods.
Assumption~\ref{ass:cost-region} implies $0<\alpha_h<1$
and $\alpha_b>0$, but does not impose an ordering between
them; moreover, $\alpha_b$ may exceed one.
Their relative ordering determines which intermediate
recovery regime arises as $\alpha$ increases.

Unless otherwise specified, the subsequent analysis is based on
Assumption~\ref{ass:cost-region}. The standing cost ordering also
guarantees a finite-horizon order-up-to structure before any explicit
envelope construction is imposed.

\begin{proposition}
\label{prop:martingale-base-stock}
For every finite $t$, $B_t(\cdot,\mu)$ is convex on $[0,U]$. Define
\begin{equation*}
\chi_t(\mu):=
\max\arg\min_{x\in[0,U]}B_t(x,\mu).
\end{equation*}
Then $x_t^*(y,\mu)=\max\{y,\chi_t(\mu)\}$ is optimal for every $(y,\mu)\in[0,U]^2$.
\end{proposition}
Here $\chi_t$ selects the largest unconstrained minimizer. In the
no-replenishment results at $\alpha=\alpha_b$, we explicitly use the
alternative optimal selection $x_t^*=y$; it need not be the largest
optimal action. Away from such stated exceptions, the displayed policies
use the largest-minimizer convention.
Thus, the recovery regimes developed next determine the explicit form and
horizon evolution of the target $\chi_t$, rather than the existence of an
order-up-to policy itself.

\section{Two-Period Envelope Geometry and Policy Structure}
\label{sec:two-period}

This section examines how automatic local recovery affects optimal replenishment decisions. We characterize the optimal policies and corresponding worst-case demand distributions for the terminal- and two-period problems. These results provide the foundation for the finite-horizon analysis developed in Section~\ref{sec:finite-horizon}.

\subsection{Terminal-period benchmark}
\label{subsec:Terminal-period solution}

The terminal-period problem is common to all four recovery regimes.
Define
\begin{equation*}
\Phi_1:=c+h-\alpha w,\quad K_1:=b+h,\quad\rho_1:=\frac{\Phi_1}{K_1},\quad\bar q_2:=\rho_1U\in(0,U),
\end{equation*}
where the strict bounds follow from Assumption~\ref{ass:cost-region}.
This threshold determines the terminal-period replenishment policy
and plays a central role in the subsequent two-period analysis.

\begin{proposition}
\label{prop:terminal-solution}
For $(y_1,\mu_1)\in[0,U]^2$, the Bellman target is
$\chi_1(\mu_1)=U\mathbf 1_{\{\mu_1\geq\bar q_2\}}$.
An optimal action and an attaining worst-case distribution are
\begin{equation*}
x_1^*(y_1,\mu_1)
=
\max\{y_1,\chi_1(\mu_1)\},
\qquad
\mathbb Q_1^*
=
\left(1-\frac{\mu_1}{U}\right)\delta_0
+
\frac{\mu_1}{U}\delta_U.
\end{equation*}
The corresponding optimal value is
\begin{equation*}
\hat V^1(y_1,\mu_1)
=
-cy_1+(c+h)U-h\mu_1
-\bigl(U-x_1^*(y_1,\mu_1)\bigr)
\left(\Phi_1-K_1\frac{\mu_1}{U}\right).
\end{equation*}
\end{proposition}

With two periods remaining, realized demand $D$ becomes the
conditional mean in the terminal period. Proposition~\ref{prop:terminal-solution}
shows that the terminal policy switches at $D=\bar q_2$,
which introduces a potential breakpoint in the continuation
cost $F_2^x(D)$.

By Lemma~\ref{lem:moment-envelope}, an attaining worst-case
distribution can be supported on at most two demand points.
For $0\leq D^-\leq\mu\leq D^+\leq U$ with $D^-<D^+$,
the corresponding expected cost defines the chord
\begin{equation*}
\ell_{D^-,D^+}^x(\mu)
:=
\frac{D^+-\mu}{D^+-D^-}F_2^x(D^-)
+
\frac{\mu-D^-}{D^+-D^-}F_2^x(D^+).
\end{equation*}
Consequently,
\begin{equation*}
\begin{aligned}
\sup_{\mathbb Q\in\mathcal M(\mu)}
\mathbb E_{\mathbb Q}[F_2^x(\widetilde D)]
=\sup_{\substack{0\leq D^-\leq\mu\leq D^+\leq U\\D^-<D^+}}
\ell_{D^-,D^+}^x(\mu).
\end{aligned}
\end{equation*}
The maximizing endpoints determine an attaining worst-case distribution.

\subsection{Two-period worst-case demand distributions}

Recall that the terminal-period order-up-to target
switches at $\bar q_2$. Since realized demand $D$
becomes the terminal reference mean, this switch
also divides the continuation cost $F_2^x(D)$
into two convex segments. Specifically, for any
fixed $x\in[0,U]$,
\begin{equation*}
F_2^x(D)\leq
\begin{cases}
\ell_{0,\bar q_2}^x(D),
& 0\leq D\leq\bar q_2,\\[1mm]
\ell_{\bar q_2,U}^x(D),
& \bar q_2\leq D\leq U.
\end{cases}
\end{equation*}
Thus, any demand realization can be replaced by
a mean-preserving mixture of the endpoints of its
respective interval without decreasing the expected
cost. It follows that an attaining worst-case
distribution can be supported on
$\{0,\bar q_2,U\}$.

It remains to determine whether $\bar q_2$ is
needed in this support. If
$F_2^x(\bar q_2)\leq\ell_{0,U}^x(\bar q_2)$,
the endpoint chord dominates $F_2^x$ throughout
$[0,U]$. Otherwise, the upper envelope consists
of the two chords meeting at $\bar q_2$.
Define $\Phi_2:=c+(1+\beta)(h-\alpha w)$ and
$K_2:=b+h+\beta\Phi_1-c$. The difference is
\begin{equation*}
\Delta_2(x):=
F_2^x(\bar q_2)-\ell_{0,U}^x(\bar q_2)
=
\begin{cases}
\displaystyle
\left(1-\frac{\bar q_2}{U}\right)
(\beta\Phi_1U-K_2x),
&x\leq\bar q_2,\\[2mm]
\displaystyle
(\beta K_1-K_2)\frac{\bar q_2}{U}(U-x),
&x\geq\bar q_2.
\end{cases}
\end{equation*}
When $\Delta_2(x)>0$, the candidate support points must include
$\bar q_2$ in addition to $0$ and $U$; otherwise, the two endpoints
$\{0,U\}$ suffice.

The sign of $\Delta_2(x)$ yields two cases. If
$K_2>\beta K_1$, $\Delta_2(x)$ changes sign at
$\bar x_2:=\beta\Phi_1U/K_2<\bar q_2$. The worst-case
distribution may involve $\bar q_2$ for $x<\bar x_2$,
whereas the endpoints $\{0,U\}$ suffice for
$x\geq\bar x_2$. If $K_2\leq\beta K_1$,
$\Delta_2(x)\geq0$ for all $x\in[0,U]$, so $\bar q_2$
remains a candidate support point regardless of $x$. The corresponding worst-case expectations are characterized
in the following lemma.

\begin{lemma}[Two-period worst-case expectation]
\label{lem:two-period-envelope}
For every $x,\mu\in[0,U]$, the worst-case expectation
of $F_2^x$ has the following structure.

If $K_2>\beta K_1$, then
\begin{equation}
\sup_{\mathbb Q\in\mathcal M(\mu)}
\mathbb E_{\mathbb Q}[F_2^x(\widetilde D)]=
\begin{cases}
\ell_{0,\bar q_2}^x(\mu),
&0\leq x<\bar x_2,\quad 0\leq\mu\leq\bar q_2,\\
\ell_{\bar q_2,U}^x(\mu),
&0\leq x<\bar x_2,\quad \bar q_2\leq\mu\leq U,\\
\ell_{0,U}^x(\mu),
&\bar x_2\leq x\leq U,\quad 0\leq\mu\leq U.
\end{cases}
\label{eq:two-period-envelope}
\end{equation}
If $K_2\leq\beta K_1$, then
\begin{equation}
\sup_{\mathbb Q\in\mathcal M(\mu)}
\mathbb E_{\mathbb Q}[F_2^x(\widetilde D)]=
\begin{cases}
\ell_{0,\bar q_2}^x(\mu),
&0\leq\mu\leq\bar q_2,\\
\ell_{\bar q_2,U}^x(\mu),
&\bar q_2\leq\mu\leq U,
\end{cases}
\qquad x\in[0,U].
\label{eq:two-period-envelope-alternative}
\end{equation}
\end{lemma}

Lemma~\ref{lem:two-period-envelope} identifies two
possible worst-case demand structures. When
$K_2>\beta K_1$, the intermediate demand level
$\bar q_2$ is relevant for $x<\bar x_2$, whereas
the endpoint distribution on $\{0,U\}$ is
worst-case for $x\geq\bar x_2$.
When $K_2\leq\beta K_1$, an attaining distribution
can retain $\bar q_2$ as a support point for any
$x$, paired with either $0$ or $U$ according to
$\mu$. The associated probabilities follow from
Lemma~\ref{lem:moment-envelope}.

Despite these differences in the worst-case
distributions, the optimal replenishment decisions
satisfy a common boundary.

\begin{corollary}[Two-period replenishment boundary]
\label{cor:penalty-boundary}
The following statements hold:
\begin{enumerate}
\item If $\alpha<\alpha_b$, positive central replenishment
is strictly beneficial in some states. In particular,
for $\mu_2=U$ and any $y_2<U$, the unique optimal
action is $x_2^*(y_2,U)=U$.
\item If $\alpha\geq\alpha_b$, no replenishment is
optimal for every $(y_2,\mu_2)\in[0,U]^2$, i.e.,
$x_2^*(y_2,\mu_2)=y_2$. This action is unique
when $\alpha>\alpha_b$.
\end{enumerate}
\end{corollary}

Corollary~\ref{cor:penalty-boundary} establishes
$\alpha_b$ as the critical recovery rate for central
replenishment in the two-period problem. Below this
threshold, replenishment can strictly reduce the
worst-case cost, whereas at or above it, no replenishment
is optimal at every state. We next characterize the
complete two-period policy for $\alpha<\alpha_b$.

\subsection{Two-Period Policy Characterization}
\label{subsec:two-period-active}

By Corollary~\ref{cor:penalty-boundary},
no replenishment is optimal when
$\alpha\geq\alpha_b$. We therefore focus on
$\alpha<\alpha_b$. Set
$\rho_2:=\Phi_2/K_2\in(0,1)$.
The optimal order-up-to target depends on
the two envelope structures identified in
Lemma~\ref{lem:two-period-envelope}.

\begin{theorem}[Two-period policy for $\alpha<\alpha_b$]
\label{thm:two-period-solution}
For every $(y_2,\mu_2)\in[0,U]^2$,
the optimal replenishment decision is $x_2^*(y_2,\mu_2)=\max\{y_2,\chi_2(\mu_2)\}$, where the order-up-to target $\chi_2$ and an attaining worst-case distribution are given by the following two cases.

\begin{enumerate}
\item If $K_2>\beta K_1$, the optimal target and an attaining worst-case distribution are
\begin{equation*}
\chi_2(\mu)=
\begin{cases}
0,
&0\leq\mu<\rho_2\bar q_2,\\
\bar x_2,
&\rho_2\bar q_2\leq\mu<\rho_2U,\\
U,
&\rho_2U\leq\mu\leq U.
\end{cases}
\end{equation*}
\begin{equation*}
\mathbb Q_2^*(\cdot\mid y_2,\mu_2)=
\begin{cases}
\displaystyle
\left(1-\frac{\mu_2}{\bar q_2}\right)\delta_0
+\frac{\mu_2}{\bar q_2}\delta_{\bar q_2},
&x_2^*<\bar x_2,\\[2mm]
\displaystyle
\left(1-\frac{\mu_2}{U}\right)\delta_0
+\frac{\mu_2}{U}\delta_U,
&x_2^*\geq\bar x_2.
\end{cases}
\end{equation*}

\item If $K_2\leq\beta K_1$, the optimal target and an attaining worst-case distribution are
\begin{equation*}
\chi_2(\mu)=
\begin{cases}
0,
&0\leq\mu<\rho_2\bar q_2,\\[1mm]
\bar q_2,
&\displaystyle
\rho_2\bar q_2\leq\mu< U-\frac{K_2-\Phi_2}{b+h-c}(U-\bar q_2),\\[2mm]
U,
&\text{otherwise}.
\end{cases}
\end{equation*}
\begin{equation*}
\mathbb Q_2^*(\cdot\mid y_2,\mu_2)=
\begin{cases}
\displaystyle
\left(1-\frac{\mu_2}{\bar q_2}\right)\delta_0
+\frac{\mu_2}{\bar q_2}\delta_{\bar q_2},
&0\leq\mu_2\leq\bar q_2,\\[2mm]
\displaystyle
\frac{U-\mu_2}{U-\bar q_2}\delta_{\bar q_2}
+\frac{\mu_2-\bar q_2}{U-\bar q_2}\delta_U,
&\bar q_2<\mu_2\leq U.
\end{cases}
\end{equation*}
\end{enumerate}
In both cases, the optimal value is
\begin{equation*}
\hat V^2(y_2,\mu_2)
=-cy_2+B_2(x_2^*(y_2,\mu_2),\mu_2).
\end{equation*}
\end{theorem}

The two cases differ in their intermediate order-up-to targets
and worst-case demand distributions. When $K_2>\beta K_1$,
the intermediate target is $\bar x_2$, and an attaining
worst-case distribution switches to the endpoint support
$\{0,U\}$ once the optimal inventory level reaches $\bar x_2$.
When $K_2\leq\beta K_1$, the intermediate target is
$\bar q_2$, which also remains a support point of an
attaining worst-case distribution. At $K_2=\beta K_1$,
we have $\bar x_2=\bar q_2$, and the upper switching
threshold in the second case equals $\rho_2U$, so the two policy expressions agree.

Together with Corollary~\ref{cor:penalty-boundary},
the preceding theorem completes the two-period
characterization under Assumption~\ref{ass:cost-region}.
The condition $K_2>\beta K_1$ distinguishes the two
possible worst-case envelope structures, whereas
$\alpha_b$ determines the replenishment boundary
independently of this distinction.
We next examine how these structures extend to
arbitrary finite horizons.

\section{Finite-Horizon Policies under the Regular Penalty Condition}
\label{sec:finite-horizon}

We extend the two-period analysis to arbitrary finite horizons.
The coefficients introduced in Section~\ref{sec:two-period} generalize to
\begin{equation}
\begin{aligned}
\Phi_t
&:=c+\sum_{k=0}^{t-1}\beta^k(h-\alpha w),
&&t\geq1,\\
K_1&:=b+h\quad\mbox{and}\quad
K_t:=b+h+\beta\Phi_{t-1}-c,
&&t\geq2.
\end{aligned}
\label{eq:finite-coefficients}
\end{equation}
These coefficients satisfy
\begin{equation*}
\begin{aligned}
&\Phi_{t+1}-\Phi_t
=\beta^t(h-\alpha w), &&t\geq1,\\
&K_t-\Phi_t
=b-c-\alpha(c-w), &&t\geq2.
\end{aligned}
\end{equation*}
The first difference determines whether $\Phi_t$
increases or decreases with the remaining horizon,
while the second corresponds to the replenishment
boundary established in Corollary~\ref{cor:penalty-boundary}.
Their signs change at $\alpha_h=h/w$ and
$\alpha_b=(b-c)/(c-w)$, respectively, yielding
the four recovery regimes in
Table~\ref{tab:recovery-regimes}.

\begin{table}[ht]
\TABLE{Recovery-rate regimes\label{tab:recovery-regimes}}
{\begin{tabular}{cccc}
\toprule
Regime I & Regime II & Regime III & Regime IV\\
\midrule
$\alpha\leq\alpha_h$
& \multirow{2}{*}{$\alpha_h<\alpha<\alpha_b$}
& \multirow{2}{*}{$\alpha_b\leq\alpha<\alpha_h$}
& \multirow{2}{*}{$\alpha\geq\max\{\alpha_h,\alpha_b\}$}\\
$\alpha<\alpha_b$ & & &\\
\bottomrule
\end{tabular}}
{}
\end{table}

The relative ordering of $\alpha_h$ and $\alpha_b$
determines which intermediate regime arises as the
recovery rate increases. Boundary cases are included as shown, and a
regime may be empty within $\alpha\in[0,1)$.

The two-period analysis in Section~\ref{sec:two-period} identifies two
possible worst-case envelope structures, distinguished
by the sign of $K_2-\beta K_1$. To develop a unified
finite-horizon analysis based on $K_2>\beta K_1$,
we focus on the following \emph{regular penalty condition}:
\begin{equation}
b-c>w-2\sqrt{hw}.
\label{eq:regular-penalty}
\end{equation}
Under Assumption~\ref{ass:cost-region}, this condition
is necessary and sufficient for $K_2>\beta K_1$ to
hold for every $\alpha\in[0,1)$. Indeed,
\begin{equation*}
K_2-\beta K_1
=(\sqrt h-\alpha\sqrt w)^2
+\alpha(b-c-w+2\sqrt{hw}).
\end{equation*}
For $t\geq3$, the corresponding inequality holds
automatically, since
\begin{equation*}
K_t-\beta K_{t-1}
=\alpha(b+h-c)
+\beta\{h+\alpha(c-w)\}>0.
\end{equation*}

\begin{figure}[ht]
    \centering
    \includegraphics[width=0.65\linewidth]{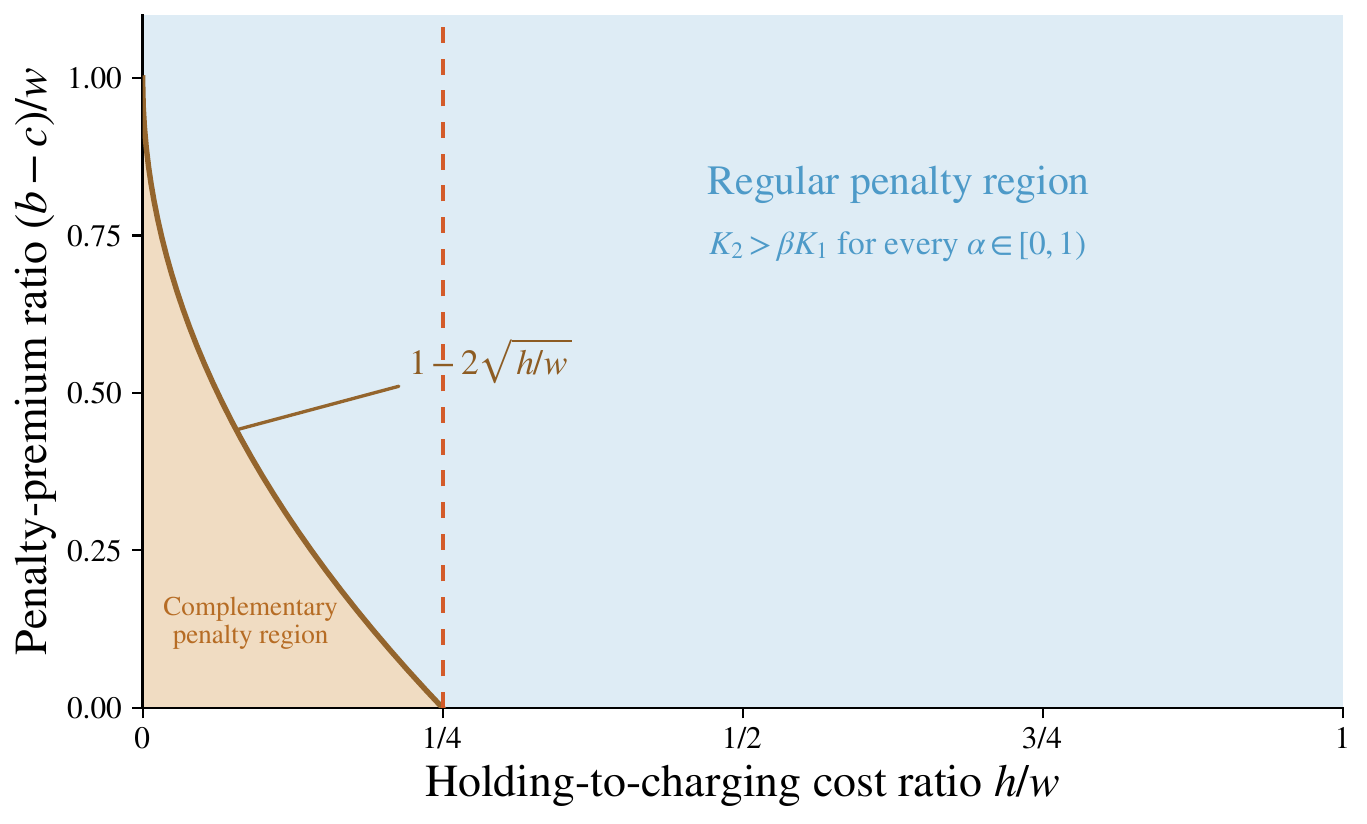}
\caption{Regular and complementary penalty regions}
\label{fig:regular-penalty-region}
\end{figure}

The regular penalty condition is economically reasonable and covers a broad range of cost configurations. It holds automatically whenever $h\geq w/4$ and, even when $h<w/4$, requires only a shortage-penalty premium $b-c$ exceeding $w-2\sqrt{hw}<w$. In battery-swapping operations, unmet demand may involve lost service revenue, customer waiting, and abandonment, making a non-negligible penalty premium plausible. Such service consequences are recognized in existing studies through stockout-probability requirements~\cite{qi2023scaling}, customer waiting and abandonment~\cite{xie2024empowering}, and explicit waiting costs~\cite{sun2019optimal}. As illustrated in Figure~\ref{fig:regular-penalty-region}, the condition excludes only a restricted portion of the normalized cost-parameter space, where relatively low holding costs coincide with a small shortage-penalty premium.

Throughout this section, Assumption~\ref{ass:cost-region}
and the regular penalty condition hold.
Our analysis is divided according to the sign of
$K_t-\Phi_t=b-c-\alpha(c-w)$, which determines
whether central replenishment can be strictly beneficial
in nonterminal periods.
Specifically, Section~\ref{subsec:active-common}
derives explicit finite-horizon policies for
$\alpha<\alpha_b$ (Regimes I and II), whereas
Section~\ref{subsec:inactive-common} establishes
the optimality of no replenishment for
$\alpha\geq\alpha_b$ (Regimes III and IV).

\subsection{Finite-horizon characterization in Regimes I and II}
\label{subsec:active-common}

We first consider $\alpha<\alpha_b$, for which
$K_t-\Phi_t=b-c-\alpha(c-w)>0$ in every nonterminal period.
Building on the two-period solution, we characterize the replenishment
policy through recursive inventory thresholds and the corresponding
worst-case demand contacts. Under the regular penalty condition, this
construction yields an explicit solution for every finite horizon.

For $t\geq2$, let $\rho_t:=\Phi_t/K_t\in(0,1)$; $\rho_1$ was defined above. For $1\leq j\leq t$,
define the demand contact points
\begin{equation}
q_{t,j}:=U\prod_{k=j}^{t-1}\rho_k,
\label{eq:demand-contact-grid}
\end{equation}
where an empty product equals one. Initialize the inventory thresholds
by $s_{1,0}=0$ and $s_{1,1}=U$. For $t\geq2$, set
$s_{t,0}=0$, $s_{t,t}=U$, and
\begin{equation}
s_{t,j}:=\frac{\Phi_{t-1}}{K_t}
\bigl(s_{t-1,j}-\alpha U\bigr),
\qquad j=1,\ldots,t-1.
\label{eq:inventory-contact-grid}
\end{equation}
In particular, $q_{2,1}=\bar q_2$ and $s_{2,1}=\bar x_2$.
Some inherited inventory thresholds may become nonpositive as automatic
recovery accumulates. Define
\begin{equation*}
a_t:=\min\{j\in\{1,\ldots,t\}:s_{t,j}>0\}.
\end{equation*}
The regular penalty condition guarantees, for all finite $t$,
\begin{equation}
\begin{gathered}
0<q_{t,1}<\cdots<q_{t,t}=U,\qquad
0<s_{t,a_t}<\cdots<s_{t,t}=U,\\
s_{t,j}<q_{t,j},\qquad a_t\leq j<t.
\end{gathered}
\label{eq:automatic-grid}
\end{equation}
Thus only the positive inventory thresholds can enter the optimal
policy. For $x<U$, put
\begin{equation*}
J_t(x):=\min\{j\in\{a_t,\ldots,t\}:x<s_{t,j}\},
\qquad J_t(U):=t.
\end{equation*}

For fixed $x$, the worst-case expectation of $F_t^x(\widetilde D)$
is obtained by piecewise-linear interpolation of $F_t^x$ at
$0,q_{t,J_t(x)},\ldots,q_{t,t}=U$.
On the action interval $x\in[\chi_t(\mu),U]$, this gives the
following explicit form of the fixed-action objective, with $j=J_t(x)$:
\begin{equation}
\begin{aligned}
B_t(x,\mu)
=\;&(c+th)U+((t-1)c-th)\mu\\
&-\sum_{k=j+1}^{t}(s_{t,k}-s_{t,k-1})
\left(\Phi_t-K_t\frac{\mu}{q_{t,k}}\right)\\
&-(s_{t,j}-x)\left(\Phi_t-K_t\frac{\mu}{q_{t,j}}\right).
\end{aligned}
\label{eq:finite-fixed-action-value}
\end{equation}

\begin{theorem}[Finite-horizon characterization in Regimes I and II]
\label{thm:finite-horizon}
For $\alpha<\alpha_b$ (Regimes I and II), the optimal order-up-to target at every finite horizon $T$ is
\begin{equation}
\chi_t(\mu)=
\begin{cases}
0,&0\leq\mu<\rho_tq_{t,a_t},\\[1mm]
s_{t,j},&\rho_tq_{t,j}\leq\mu<\rho_tq_{t,j+1},
\quad j=a_t,\ldots,t-1,\\[1mm]
U,&\rho_tU\leq\mu\leq U,
\end{cases}
\label{eq:finite-target}
\end{equation}
For every $t\in[T]$ and $(y_t,\mu_t)\in[0,U]^2$, an optimal action
and the corresponding value are
\begin{equation*}
\begin{aligned}
x_t^*(y_t,\mu_t)&=\max\{y_t,\chi_t(\mu_t)\},\\
\hat V^t(y_t,\mu_t)
&=-cy_t+B_t(x_t^*(y_t,\mu_t),\mu_t).
\end{aligned}
\end{equation*}
Moreover, with $q_t^*:=q_{t,J_t(x_t^*)}$, an attaining
worst-case conditional distribution is
\begin{equation*}
\mathbb Q_t^*(\,\cdot\mid y_t,\mu_t)
=\left(1-\frac{\mu_t}{q_t^*}\right)\delta_0
+\frac{\mu_t}{q_t^*}\delta_{q_t^*},
\end{equation*}
where $q_t^*\geq\max\{x_t^*,\mu_t\}$.
\end{theorem}

The optimal post-replenishment level is nondecreasing in both current
inventory and the conditional mean. Its possible positive targets are
given by the inventory thresholds, while the corresponding demand
contacts determine the positive support of the worst-case conditional
laws. Every positive worst-case realization exhausts the ready
inventory. Since zero is absorbing under martingale demand, the
conditional laws can be concatenated into an attaining joint
martingale distribution supported on at most $T+1$ complete demand
paths. The Online Appendix establishes both the recursive grid
ordering and the concave-envelope property at arbitrary finite horizons.

For strictly positive recovery, the number of positive interior inventory
thresholds also stabilizes as the planning horizon increases. Let
$N_t:=t-a_t$ and, for $0<\alpha<\alpha_b$, define
\[
\Phi_\infty=c+\frac{h-\alpha w}{\alpha},
\qquad
\lambda_\infty=
\frac{\Phi_\infty}
{\Phi_\infty+b-c-\alpha(c-w)},
\]
and
\[
\zeta:=
\frac{\log\!\left(\alpha/(1-\beta\lambda_\infty)\right)}
{\log\lambda_\infty}.
\]

\begin{proposition}[Stabilization of inventory thresholds]
\label{prop:threshold-stabilization}
For $0<\alpha<\alpha_b$ (Regimes I and II), the number of positive interior
thresholds $N_t$ stabilizes after finitely many periods,
with
\[
N_\infty=
\begin{cases}
\lceil\zeta\rceil,
&\alpha\leq\alpha_h,\\[1mm]
1+\lfloor\zeta\rfloor,
&\alpha_h<\alpha<\alpha_b.
\end{cases}
\]
For $\alpha\leq\alpha_h$,
$N_{t+1}\in\{N_t,N_t+1\}$, whereas in Regime II
the sequence $\{N_t\}$ need not be monotone.
\end{proposition}

Thus, although the planning horizon increases, the number
of positive interior replenishment targets eventually
remains constant. The two regimes differ in how this
stabilization occurs: the number of thresholds is
nondecreasing when $0<\alpha\leq\alpha_h$, but may
decrease before stabilizing when $\alpha>\alpha_h$. At $\alpha=0$,
all interior thresholds remain positive and $N_t=t-1$, so the count
does not stabilize.

\subsection{Finite-horizon policy in Regimes III and IV}
\label{subsec:inactive-common}

We next consider $\alpha\geq\alpha_b$. Corollary~\ref{cor:penalty-boundary}
establishes the optimality of no replenishment in the two-period
problem. The following theorem extends this result to arbitrary
finite horizons.

\begin{theorem}[Finite-horizon policy in Regimes III and IV]
\label{thm:recovery-dominant-solution}
For $\alpha\geq\alpha_b$ (Regimes III and IV), every finite horizon $T$, and
$(y,\mu)\in[0,U]^2$, an optimal replenishment policy satisfies
\begin{equation*}
x_t^*(y,\mu)=y,\qquad t=2,\ldots,T.
\end{equation*}
The nonterminal action is unique when $\alpha>\alpha_b$.
In the terminal period,
\begin{equation*}
x_1^*(y,\mu)\in
\begin{cases}
\{y\}, & \mu<\rho_1U,\\
[y,U], & \mu=\rho_1U,\\
\{U\}, & \mu>\rho_1U.
\end{cases}
\end{equation*}
\end{theorem}

Unlike Regimes I and II, the optimal nonterminal policy is
independent of both the conditional mean and the remaining
horizon. This distinction follows from the replenishment
margin $b-c-\alpha(c-w)\leq0$: the benefit of avoiding
shortages no longer outweighs the cost of central supply
and the associated reduction in automatic recovery.
Although Regimes III and IV differ in the sign of
$h-\alpha w$, both admit the same no-replenishment policy.
The terminal decision retains the single-period
order-up-to structure, since no further recovery benefit
can be realized beyond the horizon.

Together, Theorems~\ref{thm:finite-horizon} and
\ref{thm:recovery-dominant-solution} characterize the
optimal replenishment policies across the four recovery
regimes under the regular penalty condition. When
$\alpha<\alpha_b$, the optimal policy follows a stepwise
order-up-to rule whose positive interior thresholds
eventually stabilize in number when $\alpha>0$. When $\alpha\geq\alpha_b$,
central replenishment is unnecessary in every nonterminal
period. Thus, $\alpha_b$ separates two fundamentally
different replenishment behaviors, while $\alpha_h$
determines the evolution of the recursive coefficients
and inventory thresholds.

The complementary case $K_2\leq\beta K_1$, which can
arise outside the regular penalty region, is addressed
in Online Appendix~\ref{app:complementary-approx}.
We develop two regime-specific approximation policies:
a conditional-mean threshold policy in Regime II and
a no-replenishment policy in Regime IV. Both recover
the exact terminal- and two-period solutions and admit
uniform bounds on their worst-case optimality gaps.
In particular, their absolute gaps remain bounded
independently of the planning horizon, while their
relative gaps converge to zero at rate $O(1/t)$.

\section{Submartingale and Supermartingale Demand}
\label{sec:directional-drift}

The martingale model specifies the next conditional mean as the latest
observed demand. We extend this restriction to allow the conditional mean
to lie above or below that reference level. The inventory dynamics and
operating costs remain unchanged, but the admissible demand distributions
now differ in their conditional means as well as their shapes.
This section retains Assumption~\ref{ass:cost-region}; the regular penalty
condition is not imposed. We continue to use the algebraic coefficients
$\Phi_t,K_t$ in \eqref{eq:finite-coefficients} and the ratio
$\rho_t=\Phi_t/K_t$.

\subsection{Directional ambiguity and order-up-to structure}

For a reference level $\mu\in[0,U]$, the one-period ambiguity sets are
\begin{equation*}
\begin{aligned}
\mathcal M^\uparrow(\mu)
&:=\bigl\{\mathbb Q\in\mathcal P_0([0,U]):
\mathbb E_{\mathbb Q}[\widetilde D]\geq\mu\bigr\},\\
\mathcal M^\downarrow(\mu)
&:=\bigl\{\mathbb Q\in\mathcal P_0([0,U]):
\mathbb E_{\mathbb Q}[\widetilde D]\leq\mu\bigr\}.
\end{aligned}
\end{equation*}
They correspond to submartingale and supermartingale demand,
respectively. As in Section~\ref{Model}, $\mu_T=\mu$ and
$\mu_t=D_{t+1}$ for $t<T$. The associated joint ambiguity sets are
\begin{equation*}
\begin{aligned}
\mathcal{MAR}^{\uparrow}(\mu)
&:=\bigl\{\mathbb P\in\mathcal P_0([0,U]^T):
\mathbb E_{\mathbb P}[\widetilde D_t\mid\mathcal F_{t+1}]
\geq\widetilde D_{t+1},\ t\in[T]\bigr\},\\
\mathcal{MAR}^{\downarrow}(\mu)
&:=\bigl\{\mathbb P\in\mathcal P_0([0,U]^T):
\mathbb E_{\mathbb P}[\widetilde D_t\mid\mathcal F_{t+1}]
\leq\widetilde D_{t+1},\ t\in[T]\bigr\},
\end{aligned}
\end{equation*}
where $\widetilde D_{T+1}\equiv\mu$. Both sets retain the rectangular
conditional structure of the martingale model.

Let $\hat V_\uparrow^t$ and $\hat V_\downarrow^t$ denote the resulting
cost-to-go functions, initialized at
$\hat V_\uparrow^0=\hat V_\downarrow^0=0$.
The two-point attainment property and Bellman recursion extend to
both directional ambiguity sets.

\begin{proposition}[Directional dynamic programming]
\label{prop:directional-dp}
For either $\diamond\in\{\uparrow,\downarrow\}$, the corresponding
value function satisfies
\begin{equation}
\begin{aligned}
\hat V_\diamond^t(y,\mu)
=\min_{x\in[y,U]}\Bigl\{
&c(x-y)+w\alpha(U-x)\\
&+\sup_{\mathbb Q\in\mathcal M^\diamond(\mu)}
\mathbb E_{\mathbb Q}\bigl[
 b(\widetilde D-x)_{+}+h(x-\widetilde D)_{+}\\
&\hspace{37mm}
+\hat V_\diamond^{t-1}(R(x,\widetilde D),\widetilde D)
\bigr]\Bigr\}.
\end{aligned}
\label{eq:directional-bellman}
\end{equation}
For each fixed action, the supremum is attained by a distribution
supported on at most two demand points.
\end{proposition}

\begin{corollary}
\label{cor:directional-monotonicity}
For every finite $t$ and $(y,\mu)\in[0,U]^2$,
\[
\hat V^t(y,\mu)\leq\hat V_\uparrow^t(y,\mu),
\qquad
\hat V^t(y,\mu)\leq\hat V_\downarrow^t(y,\mu).
\]
Furthermore, $\hat V_\uparrow^t(y,\mu)$ is nonincreasing in $\mu$,
whereas $\hat V_\downarrow^t(y,\mu)$ is nondecreasing in $\mu$.
\end{corollary}

The opposite monotonicities reflect the admissible mean intervals:
increasing $\mu$ restricts $\mathcal M^\uparrow(\mu)$ but enlarges
$\mathcal M^\downarrow(\mu)$. In particular,
$\mathcal M^\uparrow(U)=\{\delta_U\}$ and
$\mathcal M^\downarrow(0)=\{\delta_0\}$.

Write $B_{\diamond,t}(x,\mu)$ for the fixed-action objective in
\eqref{eq:directional-bellman} after separating $-cy$, analogously
to $B_t$ in Section~\ref{Model}. The convexity argument underlying
Proposition~\ref{prop:martingale-base-stock} applies to both
conditional mean restrictions.

\begin{proposition}
\label{prop:directional-base-stock}
For either $\diamond\in\{\uparrow,\downarrow\}$ and every finite $t$,
$B_{\diamond,t}(\cdot,\mu)$ is convex on $[0,U]$. Define
\[
\chi_t^\diamond(\mu)
:=\max\arg\min_{x\in[0,U]}B_{\diamond,t}(x,\mu).
\]
Then an optimal replenishment decision is
\[
x_{\diamond,t}^*(y,\mu)
=\max\{y,\chi_t^\diamond(\mu)\}.
\]
\end{proposition}

The directional models therefore preserve the order-up-to structure,
but their explicit policies differ markedly.

\subsection{Submartingale demand: upward drift}
\label{subsec:submartingale}

Under upward ambiguity, $U$ is absorbing and $\delta_U$ belongs
to $\mathcal M^\uparrow(\mu)$ at every reference level.
When $\alpha<\alpha_b$, this structure leads to an explicit policy
and a finite horizon beyond which replenishment to capacity is optimal
at every state.

Recall $\rho_t=\Phi_t/K_t$. Define
\begin{equation*}
\omega_t:=1-\frac{\bigl(th-(t-1)c\bigr)_{+}}{K_t},
\qquad
C_t^\uparrow:=tcU+\bigl(th-(t-1)c\bigr)_{+}(1-\rho_t)U.
\end{equation*}

\begin{theorem}[Finite-horizon policy under submartingale demand]
\label{thm:submartingale-solution}

For $\alpha<\alpha_b$, the largest optimal order-up-to target
at every finite horizon $t$ is
\begin{equation}
\chi_t^\uparrow(\mu)=
\begin{cases}
\omega_tU, & 0\leq\mu<\rho_tU,\\[1mm]
U, & \rho_tU\leq\mu\leq U.
\end{cases}
\label{eq:submartingale-target}
\end{equation}
An optimal action is
\[
x_{\uparrow,t}^*(y,\mu)
=\max\{y,\chi_t^\uparrow(\mu)\},
\]
and the corresponding value is
\begin{equation*}
\begin{aligned}
\hat V_\uparrow^t(y,\mu)
={}&-cy+C_t^\uparrow
+(\Phi_t-K_t)
\bigl(x_{\uparrow,t}^*(y,\mu)-\omega_tU\bigr)\\
&+K_t\left(1-\frac{\mu}{U}\right)
\bigl(x_{\uparrow,t}^*(y,\mu)-\omega_tU\bigr)_{+}.
\end{aligned}
\end{equation*}

For $\alpha\geq\alpha_b$, no replenishment is optimal
in the two-period problem:
\[
x_{\uparrow,2}^*(y,\mu)=y,
\qquad (y,\mu)\in[0,U]^2.
\]
This action is unique when $\alpha>\alpha_b$; at equality we select
$x_{\uparrow,2}^*=y$, which need not be the largest optimal action.
\end{theorem}

The lower target $\omega_tU$ reaches capacity precisely when
$(t-1)c\geq th$. Since $c>h$, let
\begin{equation*}
t_\uparrow:=\left\lceil\frac{c}{c-h}\right\rceil.
\end{equation*}
The two branches of \eqref{eq:submartingale-target} then coincide
for every $t\geq t_\uparrow$.

\begin{corollary}
\label{cor:submartingale-capacity}
Under $\alpha<\alpha_b$, every $t\geq t_\uparrow$ satisfies
\[
x_{\uparrow,t}^*(y,\mu)=U,
\qquad
\hat V_\uparrow^t(y,\mu)=tcU-cy,
\]
for all $(y,\mu)\in[0,U]^2$.
\end{corollary}

In particular, $t_\uparrow=2$ whenever $c\geq2h$.
The capacity result does not extend to $\alpha\geq\alpha_b$:
the two-period decision instead admits a no-replenishment policy.

Thus the two-period submartingale policy shares the replenishment
boundary of Corollary~\ref{cor:penalty-boundary}. At longer horizons,
Proposition~\ref{prop:directional-base-stock} guarantees an order-up-to
policy, but does not establish no replenishment throughout the
region $\alpha\geq\alpha_b$.

\subsection{Supermartingale demand: downward drift}
\label{subsec:supermartingale}

Under downward ambiguity, zero is absorbing because
$\mathcal M^\downarrow(0)=\{\delta_0\}$. At this reference level,
additional replenishment is unnecessary; hence capacity
replenishment cannot be optimal at every state for any finite horizon.

The terminal target illustrates the difference from upward ambiguity.
Let
\[
\alpha_d:=\frac{b}{b+h}\in(0,1).
\]
Its largest optimal order-up-to target is
\begin{equation*}
\chi_1^\downarrow(\mu)=
\begin{cases}
0, & 0\leq\mu<\rho_1U,\\[1mm]
\alpha_dU, & \rho_1U\leq\mu\leq U.
\end{cases}
\end{equation*}
Even at high reference demand, the terminal target remains below
capacity because the worst-case distribution may select a lower mean.

A different effect arises when recovery is sufficiently strong.
If $\alpha>\alpha_d$, leaving the ready inventory at zero produces
$\alpha U>\alpha_dU$ units through automatic recovery, exceeding
the terminal target. Central replenishment can then reduce subsequent
recovery by shrinking the depleted-battery pool. This mechanism may
induce a positive nonterminal target even when $\alpha>\alpha_b$.

For $\alpha>\max\{\alpha_b,\alpha_d\}$, set
\begin{equation*}
\begin{aligned}
z_2^\downarrow
&:=U\left(1-\frac{\alpha_d}{\alpha}\right),\\
\gamma_0&:=c+(2-\alpha)(h-\alpha w),\\
\gamma_1&:=c-b-\alpha h-\alpha\beta w,\\
\mu_2^\dagger&:=U\frac{\gamma_0}{\gamma_0-\gamma_1}.
\end{aligned}
\end{equation*}
The standing cost ordering gives
$\gamma_0>0>\gamma_1$, and hence $\mu_2^\dagger\in(0,U)$.

\begin{proposition}[Two-period recovery suppression under downward drift]
\label{prop:supermartingale-recovery-suppression}
Under $\alpha>\max\{\alpha_b,\alpha_d\}$, the largest optimal
two-period target is
\begin{equation*}
\chi_2^\downarrow(\mu)=
\begin{cases}
0, & 0\leq\mu<\mu_2^\dagger,\\[1mm]
z_2^\downarrow, & \mu_2^\dagger\leq\mu\leq U.
\end{cases}
\end{equation*}
Consequently,
\[
x_{\downarrow,2}^*(y,\mu)
=\max\{y,\chi_2^\downarrow(\mu)\}.
\]
\end{proposition}

Thus $\alpha_b$ does not by itself separate positive and zero
replenishment under downward ambiguity: a sufficiently high recovery
rate can make replenishment valuable by reducing future recovery.

We next obtain an explicit finite-horizon policy for
$\alpha<\alpha_b$ and
$\alpha\leq\min\{\alpha_h,\alpha_d\}$.
The construction retains the demand contacts $q_{t,j}$ from
\eqref{eq:demand-contact-grid} but uses a different recursion for
inventory targets. Initialize
$s_{1,0}^\downarrow=0$ and
$s_{1,1}^\downarrow=\alpha_dU$. For $t\geq2$, set
$s_{t,0}^\downarrow=0$ and
\begin{equation}
\begin{aligned}
s_{t,j}^\downarrow
&:=\min\left\{U,
\left[\frac{\Phi_{t-1}}{K_t}
\bigl(s_{t-1,j}^\downarrow-\alpha U\bigr)\right]_{+}\right\},
\quad j=1,\ldots,t-1,\\
s_{t,t}^\downarrow
&:=\min\left\{U,
\left[\frac{\Phi_{t-1}}{K_t}
\bigl(s_{t-1,t-1}^\downarrow-\alpha U\bigr)
+\frac{b}{K_t}U\right]_{+}\right\}.
\end{aligned}
\label{eq:downward-inventory-grid}
\end{equation}

With $\hat V_\downarrow^0\equiv0$, the directional fixed-action objective
on $x\in[\chi_t^\downarrow(\mu),U]$ has the explicit form
\begin{equation}
\begin{aligned}
B_{\downarrow,t}(x,\mu)
={}&(c+th)U
+\left(\frac{\hat V_\downarrow^{t-1}(0,U)}{U}-th\right)_{+}\mu\\
&-\Phi_t\bigl(U-\max\{x,s_{t,t}^\downarrow\}\bigr)\\
&-\sum_{j=1}^{t}
\bigl(s_{t,j}^\downarrow
-\max\{x,s_{t,j-1}^\downarrow\}\bigr)_{+}
\left(\Phi_t-K_t\frac{\mu}{q_{t,j}}\right).
\end{aligned}
\label{eq:downward-value}
\end{equation}

\begin{theorem}[Finite-horizon policy under downward drift]
\label{thm:supermartingale-solution}
Under $\alpha<\alpha_b$ and
$\alpha\leq\min\{\alpha_h,\alpha_d\}$, the largest optimal
order-up-to target at every finite horizon $t$ is
\begin{equation}
\chi_t^\downarrow(\mu)=
\begin{cases}
0,
&0\leq\mu<\rho_tq_{t,1},\\[1mm]
s_{t,j}^\downarrow,
&\rho_tq_{t,j}\leq\mu<\rho_tq_{t,j+1},
\quad j=1,\ldots,t-1,\\[1mm]
s_{t,t}^\downarrow,
&\rho_tU\leq\mu\leq U.
\end{cases}
\label{eq:downward-target}
\end{equation}
An optimal action and its value are
\begin{equation*}
\begin{aligned}
x_{\downarrow,t}^*(y,\mu)
&=\max\{y,\chi_t^\downarrow(\mu)\},\\
\hat V_\downarrow^t(y,\mu)
&=-cy+B_{\downarrow,t}\bigl(x_{\downarrow,t}^*(y,\mu),\mu\bigr).
\end{aligned}
\end{equation*}
\end{theorem}

Some inherited targets may be zero, so the number of positive
interior targets can be smaller than the number of recursive
candidates. Define
\[
a_t^\downarrow
:=\min\{j\in\{1,\ldots,t\}:s_{t,j}^\downarrow>0\},
\qquad
N_t^\downarrow:=t-a_t^\downarrow.
\]
The number of distinct positive target levels is
$N_t^\downarrow+1$ and has the same limiting count as in the
martingale model.

\begin{corollary}[Stabilization of supermartingale targets]
\label{cor:supermartingale-stabilization}
Under the conditions of Theorem~\ref{thm:supermartingale-solution}
and $\alpha>0$, the sequence $\{N_t^\downarrow\}$ is
nondecreasing, with
\[
N_{t+1}^\downarrow\in\{N_t^\downarrow,N_t^\downarrow+1\}.
\]
It stabilizes after a finite horizon at the value $N_\infty$
in Proposition~\ref{prop:threshold-stabilization}. Specifically,
\[
T_\downarrow^{\rm stab}
=\min\bigl\{t\geq N_\infty+1:
 s_{t,t-N_\infty}^\downarrow>0\bigr\},
\qquad
N_t^\downarrow=N_\infty\quad(t\geq T_\downarrow^{\rm stab}).
\]
\end{corollary}

Thus the positive supermartingale target levels eventually
stabilize in number, whereas upward ambiguity leads to uniform
capacity replenishment once $t\geq t_\uparrow$. The two directional
restrictions preserve the order-up-to structure but generate
distinct effects of automatic recovery on replenishment decisions.

\section{Numerical Study}
\label{Numerical}
We evaluate M-DRO against five benchmark policies using
conditional Beta demand processes with varying dispersion
and departures from the martingale mean assumption.
We first compare out-of-sample costs and service levels
across demand conditions and historical sample sizes,
then assess robustness to changes in demand distributions
between training and deployment.
Finally, we examine how central replenishment costs,
planning horizons, and automatic recovery rates affect
operating performance and replenishment decisions.

\subsection{Experimental Setting and Benchmark Policies}
\label{sec:history-setting}

We consider a baseline station with capacity $U=25$,
a planning horizon of $T=24$, initial inventory $U/2$,
and a latest observed demand of 10, with no terminal
salvage value.  Motivated by the battery-swapping setting studied in \cite{qi2023scaling}, we generate heterogeneous operating instances by independently sampling the cost and recovery parameters from uniform distributions over
\[
h\in[0.5,20],\quad w\in[30,50],\quad c\in[40,110],\quad
b\in[100,120],\quad \alpha\in[0.02,0.8].
\]
We retain samples satisfying Assumption~\ref{ass:cost-region} and the regular penalty condition~\eqref{eq:regular-penalty}, ensuring that the analytical policies apply to all generated instances. These instances span all four recovery regimes.

We generate temporally dependent demand using a conditional
Beta distribution. For each replication, we independently draw
$p_{\mathrm B}\sim\operatorname{Unif}(0,1)$, which remains fixed
throughout the replication and is independent of the operating
cost parameters. 
Conditional on $p_{\mathrm B}$ and the latest observed demand $m$,
the next-period demand is generated according to
\begin{equation*}
\begin{aligned}
q_\delta(m)
&=(1-\delta)\frac{m}{U}+\delta p_{\mathrm B},\\
\widetilde D^{\mathrm{cont}}_{j+1}
\mid(\widetilde D_j=m,p_{\mathrm B})
&\sim U\operatorname{Beta}\bigl(
\kappa q_\delta(m),
\kappa[1-q_\delta(m)]
\bigr),
\end{aligned}
\end{equation*}
where $j$ indexes periods in chronological order. The continuous draw
is rounded randomly to adjacent integers, preserving its conditional mean;
$\widetilde D_j$ denotes the resulting integer demand. Thus $j$ increases
as the remaining-horizon index $t$ decreases. Details appear in
Online Appendix~\ref{app:history-study}. The parameter $\delta$ controls mean reversion toward
$p_{\mathrm B}U$, since
$\mathbb E[\widetilde D_{j+1}\mid\widetilde D_j=m,p_{\mathrm B}]-m
=\delta(p_{\mathrm B}U-m)$, whereas $\kappa>0$ determines conditional demand dispersion. The
conditional variance of the continuous draw is
$U^2q_\delta(m)[1-q_\delta(m)]/(\kappa+1)$; random rounding adds
a conditional variance term of at most $1/4$. The mean-reversion
parameter satisfies $0\leq\delta\leq1$. We consider
$\delta\in\{0.075,0.15,0.30\}$ and
$\kappa\in\{0.125,0.25,0.5,0.75,1,1.5,2,3,4,6,8\}$,
covering different degrees of mean reversion and
conditional dispersion, including extreme-dispersion
stress tests at small values of $\kappa$.

We compare M-DRO with five benchmark policies that differ
in their use of demand information and replenishment
decision rules. All policies observe the same current
inventory and latest demand, but differ in whether and
how they use historical observations, as summarized
in Table~\ref{tab:benchmark-policies}.

\begin{table}[ht]
\TABLE{Policies compared in the numerical study
\label{tab:benchmark-policies}}
{\small\begin{tabular}{lcl}
\toprule
Policy & Historical fitting & Decision rule\\
\midrule
M-DRO & No & Analytical DRO policy\\
Rolling CE & No & Rolling optimization with latest-demand forecast\\
\midrule
IID DP & Yes & DP with pooled empirical demand\\
Conditional DP & Yes & DP with empirical conditional transitions\\
Fixed Base-stock & Yes & Data-fitted constant order-up-to level\\
Myopic Newsvendor & Yes & Newsvendor with kernel-estimated quantiles\\
\bottomrule
\end{tabular}}
{}
\end{table}

Specifically, rolling CE uses the latest observed demand as a point
forecast, allowing us to isolate the value of accounting
for demand uncertainty without historical fitting.
IID DP and Conditional DP assess the benefits of
estimating unconditional and state-dependent demand
distributions, respectively. Fixed Base-stock and
Myopic Newsvendor provide simpler alternatives that
do not optimize over the full planning horizon.
Together, these benchmarks distinguish the effects
of distribution estimation, temporal dependence,
and dynamic optimization.

All historical fitting procedures are fixed before
out-of-sample evaluation, and no policy has access
to the true demand-generating parameters.

We generate 80 independent operating instances, reuse them across
demand settings, and evaluate all policies on
250 test paths per instance. The data-driven benchmarks
are trained using $N\in\{24,48,72,96\}$ historical
demand transitions collected from independent
24-period episodes. Training and test data are
independent, while all policies are evaluated on
identical test paths within each instance. Percentage cost savings are calculated relative
to each benchmark's mean cost.
The DP benchmarks use a 0.125-unit inventory grid,
with interpolation for off-grid states and an
explicit no-replenishment option.

\subsection{Out-of-Sample Performance}
\label{sec:out-of-sample}

\begin{figure}[ht]
\centering
\includegraphics[width=0.85\linewidth]{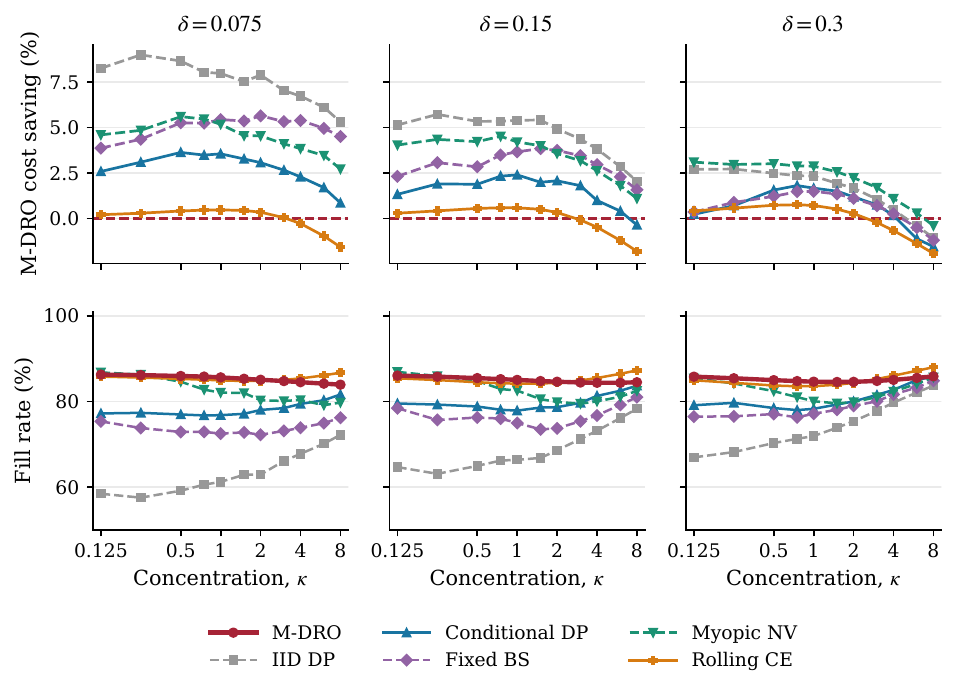}
\caption{Out-of-sample cost savings and fill rates under persistent Beta demand. ($N=48$, $T=24$)}
\label{fig:beta-oos}
\end{figure}

We first examine out-of-sample performance when training
and test data follow the same demand distribution.
Figure~\ref{fig:beta-oos} compares the six policies in
terms of operating costs and fill rates.
Across all tested values of $\delta$, M-DRO achieves
the lowest sample mean cost when $\kappa\leq2$,
with this advantage extending to $\kappa=3$ at
$\delta=0.075$. Its cost savings are generally larger
against the fitted benchmarks than against Rolling CE.
This pattern suggests that when conditional demand
is highly dispersed, estimating its distribution from
limited historical observations may offer less value
than making decisions directly from the latest demand
information.

For instance, at $\delta=0.075$ and $\kappa=0.25$,
M-DRO reduces mean cost by 3.10\%--9.01\% relative
to the four fitted benchmarks, compared with 0.29\%
against Rolling CE ($N=48$). The paired 95\%
confidence intervals for savings against Conditional
DP and Rolling CE are $[1.32\%,4.88\%]$ and
$[0.20\%,0.38\%]$, respectively.
The relatively small difference from Rolling CE
suggests that adapting replenishment decisions to
recent demand observations can already provide
substantial operational value without distribution
estimation. As $\kappa$ increases, however, M-DRO's
cost advantage diminishes and may reverse; at
$\delta=0.075$ and $\kappa=3$, its difference from
Rolling CE is no longer statistically significant.
Moreover, M-DRO does not consistently achieve
higher fill rates, indicating that cost efficiency
and service performance must be assessed separately.

\begin{table}[ht]
\centering
\TABLE{Effect of historical sample size on operating cost
\label{tab:beta-oos}}
{\small\setlength{\tabcolsep}{6pt}
\begin{tabular}{@{}crrrrrr@{}}
\toprule
& \multicolumn{2}{c}{No historical fitting}
& \multicolumn{4}{c}{Conditional DP}\\
\cmidrule(lr){2-3}\cmidrule(l){4-7}
$\kappa$ & M-DRO & Rolling CE
& $N=24$ & $N=48$ & $N=72$ & $N=96$\\
\midrule
0.125 & \textbf{910.6} & 913.2 & 946.6 & 922.9 & 915.3 & 915.1\\
0.25  & \textbf{905.0} & 908.8 & 947.2 & 922.5 & 912.4 & 912.3\\
0.5   & \textbf{894.4} & 899.4 & 928.9 & 911.6 & 904.5 & 902.5\\
0.75  & \textbf{884.5} & 889.7 & 922.6 & 905.5 & 895.4 & 893.3\\
1     & \textbf{875.7} & 880.8 & 913.3 & 897.2 & 885.8 & 884.8\\
1.5   & \textbf{861.7} & 866.0 & 897.2 & 879.2 & 871.1 & 869.5\\
2     & \textbf{850.0} & 852.8 & 884.7 & 868.0 & 859.3 & 856.8\\
3     & 831.9 & \textbf{831.2} & 856.7 & 847.2 & 838.7 & 835.7\\
8     & 789.2 & \textbf{775.2} & 795.1 & 786.4 & 780.2 & 777.5\\
\bottomrule
\end{tabular}}
{$\delta=0.15$; mean cost per period. Conditional DP: the best-performing fitting benchmark}
\end{table}

Table~\ref{tab:beta-oos} compares M-DRO and Rolling CE,
which require no historical fitting, with Conditional DP,
the best-performing fitted benchmark in these settings,
trained on different numbers of historical demand
observations. Increasing the sample size consistently
reduces Conditional DP's operating cost, but its
relative performance depends strongly on demand
dispersion. For all $\kappa\leq2$, M-DRO retains the
lowest mean cost even when Conditional DP uses
96 observations. At $\kappa=0.125$, for example,
increasing $N$ from 24 to 96 reduces Conditional DP's
cost from 946.6 to 915.1, still above M-DRO's 910.6.
When demand becomes more concentrated, however,
this advantage diminishes. At $\kappa=8$, Conditional
DP outperforms M-DRO with only 48 observations,
although Rolling CE remains the least costly policy.
Rolling CE also slightly outperforms M-DRO at
$\kappa=3$.

These findings suggest that the value of historical
demand estimation depends on demand variability
as well as sample size. Under highly dispersed demand,
even substantial increases in historical observations
do not necessarily make fitted policies preferable
to M-DRO. When demand is more concentrated,
distribution estimation becomes more effective,
although Rolling CE may still provide a better
alternative without historical fitting.
For inventory managers, this highlights the
importance of selecting replenishment policies
according to the demand environment rather than
simply increasing the amount of historical data.

\subsection{Robustness Testing}
\label{sec:beta-robustness}

We further examine policy robustness to demand shifts across heterogeneous
historical environments. In each of 80 independent replications,
we draw $\delta\sim\mathrm{Unif}(0.075,0.30)$ and
$\kappa\sim\mathrm{Unif}(0.5,2)$, independently of the operating
parameters, and train the four fitted benchmarks using $N=48$
historical transitions. The fitted policies are then held fixed.
We consider high, medium, and low dispersion settings, corresponding
to $\kappa/2$, $2\kappa$, and $4\kappa$, respectively; changes in
mean reversion ($\delta+0.10$ and $\delta-0.05$); and joint shifts
combining $\delta+0.10$ with $\kappa/2$ or $2\kappa$.
We also include the limiting case $\delta=1$, in which demand is
conditionally independent of its latest realization.
Within each replication, operating parameters and training histories
are fixed across settings, and all policies face the same test paths
within each setting.

\begin{table}[ht]
\centering
\TABLE{Operating costs under demand distribution shifts\label{tab:beta-robustness}}
{\small\setlength{\tabcolsep}{4pt}\begin{tabular}{@{}lrrrrrr@{}}
\toprule
Demand shift & M-DRO & Rolling CE & IID DP & Cond. DP & Fixed BS & Myopic NV\\
\midrule
No shift & \textbf{872.9} & 877.3 & 908.6 & 891.2 & 896.2 & 901.6\\
High dispersion & \textbf{892.5} & 897.5 & 933.1 & 909.3 & 917.2 & 925.3\\
Medium dispersion & \textbf{846.7} & 847.4 & 876.0 & 867.1 & 868.6 & 870.3\\
Low dispersion & 818.7 & \textbf{812.3} & 840.7 & 841.1 & 840.3 & 837.7\\
Higher mean reversion & \textbf{873.6} & 878.3 & 896.5 & 889.6 & 889.0 & 899.9\\
Lower mean reversion & \textbf{870.8} & 874.9 & 915.0 & 890.8 & 899.2 & 901.7\\
High dispersion and reversion & \textbf{894.3} & 900.1 & 922.3 & 908.5 & 911.0 & 924.8\\
Medium dispersion and reversion & \textbf{847.0} & 847.6 & 864.1 & 865.2 & 862.1 & 868.2\\
\midrule
IID demand (stress test) & 869.8 & 876.4 & \textbf{848.9} & 876.4 & 858.6 & 893.4\\
\bottomrule\end{tabular}}
{$N=48$, $T=24$; mean cost per period across 80 replications.}
\end{table}

Table~\ref{tab:beta-robustness} shows that M-DRO maintains the
lowest sample mean cost under increased dispersion and the tested
moderate changes in mean reversion, including both joint shifts.
Its advantage over Rolling CE nevertheless narrows as demand becomes
more concentrated: when $\kappa$ doubles, their mean costs are
846.7 and 847.4, respectively, with a similarly small difference
under the corresponding joint shift. These narrow margins should
not be interpreted as strong evidence of superiority. When $\kappa$
quadruples, Rolling CE instead performs best, costing 812.3 versus
818.7 for M-DRO. A different reversal occurs at $\delta=1$:
IID DP costs 848.9 versus 869.8 for M-DRO and 876.4 for Conditional
DP, although both fitted DPs use historically dependent demand.

These patterns distinguish robustness to demand variability from
robustness to changes in temporal dependence. When recent demand
remains informative, M-DRO avoids fitting a potentially unstable
conditional distribution and is competitive under dispersion shifts.
As demand concentrates, Rolling CE benefits from more predictable
outcomes. When lag dependence disappears, even a previously fitted
IID model can be preferable to policies that condition heavily on the
latest observation. For managers, the central question is therefore
not merely whether demand has changed, but whether recent demand
remains a reliable signal for replenishment decisions.

\subsection{Sensitivity Analysis and Managerial Insights}
\label{sec:beta-sensitivity}

\begin{figure}[ht]
\centering
\includegraphics[width=0.85\linewidth]{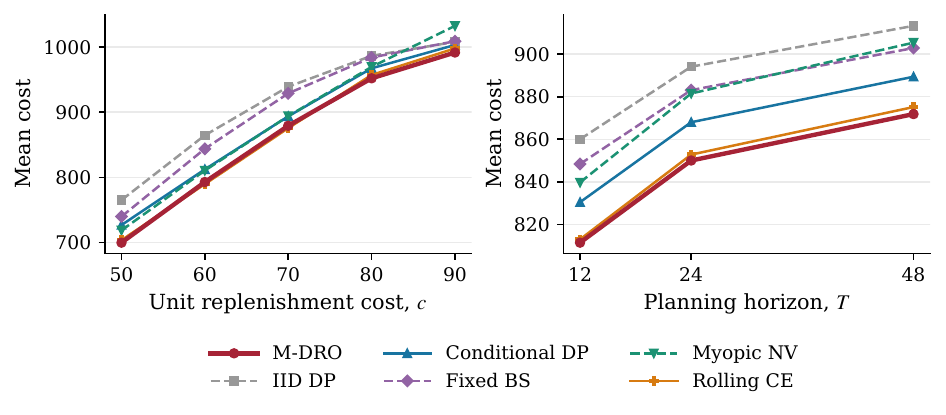}
\caption{Sensitivity to central replenishment cost and planning horizon.}
\label{fig:beta-operating}
\end{figure}

We further examine how central replenishment costs,
planning horizons, and automatic recovery rates affect
operating performance. All experiments use $N=48$,
$\delta=0.15$, and matched training and deployment
distributions. Figure~\ref{fig:beta-operating} varies
$c\in\{50,60,70,80,90\}$ and $T\in\{12,24,48\}$
at $\kappa=2$. For the replenishment-cost experiment,
we use 80 common operating instances satisfying the
cost assumptions at every tested $c$, while holding
other parameters and demand realizations fixed.

Figure~\ref{fig:beta-operating} shows that higher
central replenishment costs increase total operating
costs and encourage greater reliance on local recovery.
This behavior is consistent with the replenishment
threshold $\alpha_b=(b-c)/(c-w)$, which decreases
with $c$. As $c$ increases from 50 to 90, M-DRO's
average central replenishment decreases from 10.28
to 3.09 units per period, while local recovery
increases from 0.78 to 4.03. However, this substitution
does not fully preserve service availability: the fill
rate declines from 97.6\% to 64.3\%. M-DRO maintains
lower mean costs than the four fitted benchmarks
throughout the tested cost range, although Rolling CE
slightly outperforms it at $c=60$ and $70$.
M-DRO also achieves the lowest mean cost across all
three planning horizons, indicating that its relative
performance remains stable as the horizon increases.

\begin{figure}[ht]
\centering
\includegraphics[width=\linewidth]{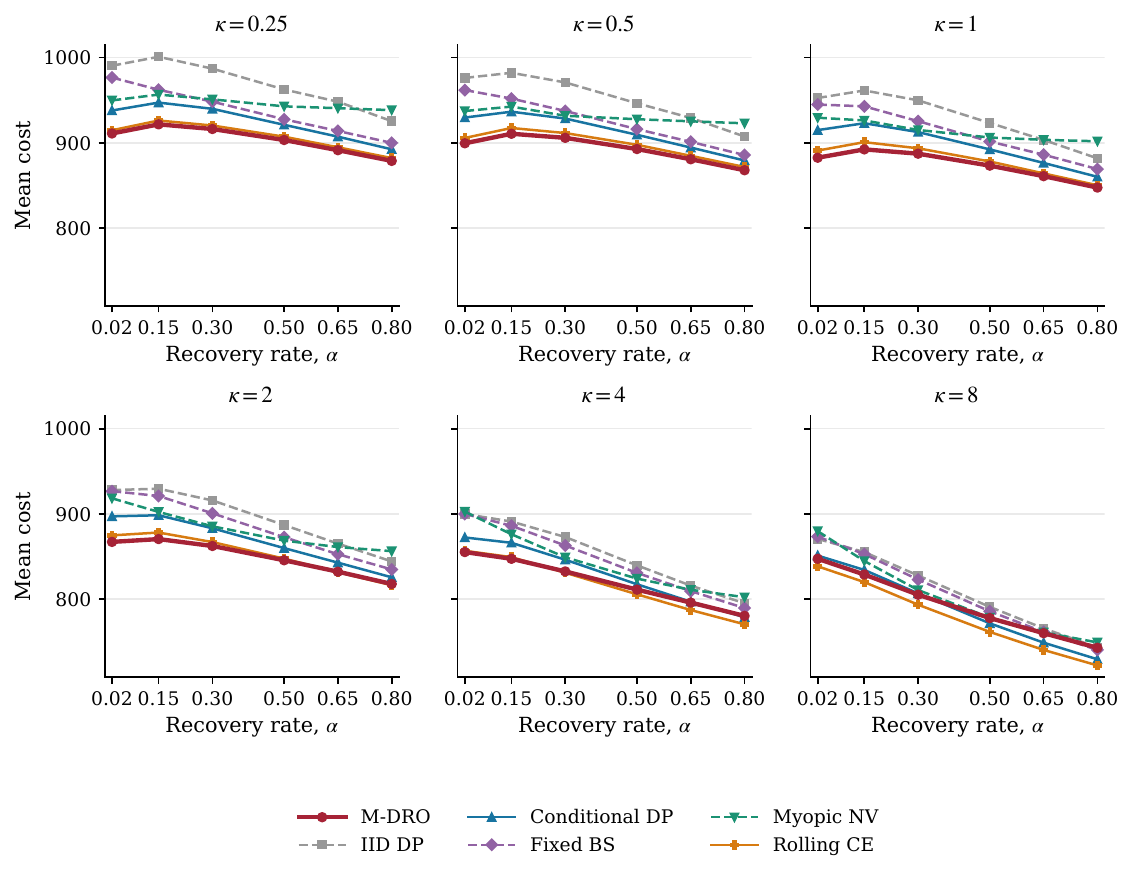}
\caption{Recovery intensity and mean operating cost.}
\label{fig:beta-recovery}
\end{figure}

Figure~\ref{fig:beta-recovery} examines the interaction
between recovery intensity and demand dispersion across
six values of $\kappa$. Operating costs generally
decline as recovery becomes stronger, although small
increases in $\alpha$ can initially raise costs under
highly dispersed demand. M-DRO achieves the lowest
mean cost throughout the tested recovery range for
$\kappa\leq1$. As demand becomes more concentrated,
its advantage over Rolling CE diminishes. At
$\kappa=2$, Rolling CE becomes the least costly policy
when $\alpha\geq0.65$, while at $\kappa=4$ this
transition occurs at $\alpha=0.30$. When $\kappa=8$,
Rolling CE consistently outperforms all alternatives.
Thus, the relative value of distributionally robust
replenishment depends jointly on recovery intensity
and demand dispersion.

Together, Figures~\ref{fig:beta-operating} and
\ref{fig:beta-recovery} reveal two important
operational trade-offs. First, substituting local
recovery for costly central replenishment can reduce central supply requirements but may compromise service availability. Second, stronger recovery does not uniformly reduce operating costs or preserve the advantage of M-DRO across demand environments. For inventory managers, effective replenishment planning therefore requires balancing recovery capability, central supply costs, and demand uncertainty rather than optimizing each consideration in isolation.

\section{Concluding Remarks}
\label{sec:conclusion}

This paper studies dynamic inventory control for battery-swapping
operations with central replenishment and automatic local recovery.
Motivated by the limited availability of historical demand data
and the difficulty of predicting temporally dependent demand,
we develop a time-consistent distributionally robust inventory
model under martingale ambiguity. We establish the finite-horizon
order-up-to structure and identify two critical recovery rates
that characterize the interplay between central replenishment
and local recovery. Under a regular penalty condition, we
derive explicit finite-horizon policies and worst-case demand
distributions, revealing when central replenishment is beneficial
and when it becomes unnecessary. Extensions to submartingale
and supermartingale demand further reveal distinct replenishment
behaviors under directional demand uncertainty.

Numerical experiments demonstrate that M-DRO is particularly
competitive when demand is highly dispersed and historical
information is limited. However, fitted or rolling-horizon
policies may perform better when demand becomes more predictable. The results
also highlight the trade-off between central replenishment,
local recovery, and service availability. Future research
could extend the model to more general demand dependence,
data-driven conditional-mean updates, and operational features
such as replenishment lead times and stochastic recovery
capacity, while preserving tractable policy structures.

\bibliographystyle{informs2014}
\bibliography{reference}
\input{appendix}
\end{document}

%% file: appendix.tex
\ECSwitch
% APPENDIX-ONLY COMPACT LAYOUT (optional): if journal requires the standard
% spacing, remove the following five lines and use the Standard version.
\linespread{1.42}\selectfont
\setlength{\abovedisplayskip}{3pt plus 1pt minus 1pt}
\setlength{\belowdisplayskip}{3pt plus 1pt minus 1pt}
\setlength{\abovedisplayshortskip}{1pt plus 1pt minus 1pt}
\setlength{\belowdisplayshortskip}{2pt plus 1pt minus 1pt}
\providecommand{\theHlemma}{}
\renewcommand{\theHlemma}{EC.\arabic{lemma}}
\providecommand{\theHproposition}{}
\renewcommand{\theHproposition}{EC.\arabic{proposition}}
\renewcommand{\theHsection}{EC.\arabic{section}}
\renewcommand{\theHsubsection}{\theHsection.\arabic{subsection}}
\renewcommand{\theHequation}{EC.\arabic{equation}}
\renewcommand{\theHtable}{EC.\arabic{table}}
\renewcommand{\theHfigure}{EC.\arabic{figure}}
\ECHead{Online Appendix: Proofs and Numerical Implementation}

Unprefixed references refer to the main manuscript. Notation follows
Section~\ref{Model}; numerical implementation appears in the final section.

\section{Dynamic-Programming Foundations}
\label{app:dp-foundations}

\proof{Proof of Proposition \ref{prop:dp-equivalence} and Corollary \ref{cor:time-consistency}.}
On the Borel demand space $[0,U]^T$, any law in $\mathcal{MAR}(\mu)$
admits conditional kernels in $\mathcal M(D_{t+1})$ almost surely;
on null histories they may be replaced by $\delta_{D_{t+1}}$
(initially $\delta_\mu$). Conversely, Ionescu--Tulcea concatenation of
Borel kernels chosen independently at each history yields precisely
$\mathcal{MAR}(\mu)$. The ambiguity set is therefore rectangular.
The correspondences $y\mapsto[y,U]$ and $m\mapsto\mathcal M(m)$
are nonempty, compact valued, and continuous in the weak topology;
for the latter, lower hemicontinuity follows by mixing a candidate law
with $\delta_0$ or $\delta_U$ when its mean is perturbed. By
finite-horizon maximum and measurable-selection theorems
(\citealt{iyengar2005robust,shapiro2016rectangular}), induction from
$\hat V^0=0$ gives continuous Bellman values, Borel minimizers
$x_t^*(y,m)$, and maximizing kernels $Q_t^*(x,m)$.
For the induced Markov policy $\pi^*$, the Bellman equation gives,
under \emph{every} admissible conditional law,
\[
\mathbb E[C(y_t,x_t^*,\widetilde D_t)+
 \hat V^{t-1}(y_{t-1},\widetilde D_t)\mid\mathcal F_{t+1}]
\leq\hat V^t(y_t,\mu_t).
\]
Iterated conditioning yields $\sup_P\mathbb E_P\sum_t C^{\pi^*}_t
\leq\hat V^T(y,\mu)$. Conversely, fix any admissible $\pi$ (for
initial inventory $y$) and compose $Q_t^*(x_t^\pi,\mu_t)$ with its
Borel history-dependent decisions. The resulting joint law belongs
to $\mathcal{MAR}(\mu)$, and Bellman optimality makes each one-step
inequality reverse. Consequently
$\sup_P\mathbb E_P\sum_t C_t^\pi\geq\hat V^T(y,\mu)$.
Taking the infimum establishes equality. For $\pi^*$ the maximizing
kernels make each step an equality, proving attainment against this
policy and optimal continuation from every state; no common
worst-case law for all policies is asserted.\Halmos\endproof

\proof{Proof of Lemma \ref{lem:moment-envelope}.}
Jensen's inequality bounds the expectation by the least concave
majorant at $\mu$; the one-dimensional concave-envelope theorem
represents this value by a contact point or a chord between two contacts
(\citealt{shapiro2012minimax}). The mean-matching Dirac or two-point
law attains the bound.\Halmos\endproof

\begin{lemma}[Inventory convexity]
\label{lem:inventory-convexity}
If $b+h\geq c$, the martingale and both directional Bellman
objectives are convex in the replenishment level at every finite horizon.
Their values are convex in inventory and satisfy
\begin{equation}
V^t(y',\mu)-V^t(y,\mu)\geq-c(y'-y),\qquad y'\geq y.
\label{eq:inventory-slope}
\end{equation}
\end{lemma}

\proof{Proof.}
Write $V$ for any of the three value functions. From
$V^t(y,\mu)+cy=\min_{x\in[y,U]}B_t(x,\mu)$, the slope bound
\eqref{eq:inventory-slope} is immediate. Inductively assume
$V^{t-1}(\cdot,D)$ is convex. The fixed-demand integrand is convex in
$x$ on either side of $x=D$; at this kink, with
$r=\alpha(U-D)$ and $v=\partial_y^+V^{t-1}(r,D)\ge-c$, its
right-minus-left derivative is $b+h+v\ge b+h-c\ge0$.
Thus $B_t(\cdot,\mu)$ is convex after expectation, supremum and affine
costs. For $y=\lambda y_1+(1-\lambda)y_2$, convex combinations of
minimizers $x_i\ge y_i$ are feasible, so
$\min_{x\ge y}B_t(x,\mu)$ is convex in $y$.
Subtracting $cy$ closes the induction; boundary cases follow by
one-sided continuity.\Halmos\endproof

\proof{Proof of Proposition \ref{prop:martingale-base-stock}.}
Convexity of $B_t(\cdot,\mu)$ implies it is nondecreasing to the right of
its largest minimizer $\chi_t(\mu)$. Its largest minimizer on $[y,U]$
is therefore $\max\{y,\chi_t(\mu)\}$.\Halmos
\endproof

\section{The Terminal and Two-Period Problems}
\label{app:two-period-proofs}

\proof{Proof of Proposition \ref{prop:terminal-solution}.}

For fixed $x$, the shortage-holding loss is convex in demand, so its least
concave majorant is the endpoint chord, attained by
$(1-\mu/U)\delta_0+(\mu/U)\delta_U$. The Bellman objective is therefore
\(-cy+w\alpha U+b\mu+x(\Phi_1-K_1\mu/U).\)
Its action slope changes sign at $\mu=\rho_1U$; minimizing over $[y,U]$
gives the stated action, and substitution gives the value.
\Halmos
\endproof

\proof{Proof of Lemma \ref{lem:two-period-envelope}.}
Write $q=\bar q_2$, $s=\bar x_2$. The terminal value shows that
$F_2^x$ is piecewise convex with breakpoints $x,q$.
For $x<q$,
\[
\ell_{0,q}^x(x)-F_2^x(x)
=\frac{K_1(b+2h-\alpha w)}{U\Phi_1}x(q-x)\ge0;
\]
$F_2^x$ is affine on $[q,U]$, while the difference between
the left and right chord slopes is
$K_2(s-x)/q$. Thus the adjacent chords are concave
for $x<s$; for $s\le x<q$, the endpoint chord dominates if
$K_2>\beta K_1$ because
$\ell_{0,U}^x(q)-F_2^x(q)=K_2(1-q/U)(x-s)\ge0$.
When $K_2\le\beta K_1$, $s\ge q$, and adjacent chords remain.
For $x\ge q$, convexity holds on both sides of $q$ (the derivative
jump at $D=x$ is $b+h-c>0$), and
\[
\ell_{0,U}^x(q)-F_2^x(q)
=(K_2-\beta K_1)\frac qU(U-x).
\]
The sign determines the corresponding least concave majorant:
all chosen chords dominate and touch $F_2^x$.
Lemma~\ref{lem:moment-envelope} gives the asserted envelopes.\Halmos\endproof

\proof{Proof of Corollary \ref{cor:penalty-boundary}.}
At $\mu=U$, demand remains $U$ and the two-period objective is
$2cU-cy+[b-c-\alpha(c-w)](U-x)$; thus $x=U$ is uniquely optimal
when $\alpha<\alpha_b$.
For $\alpha\ge\alpha_b$, let
$\eta=\Phi_2-K_2=(c-w)(\alpha-\alpha_b)\ge0$.
All one-sided slopes obtained from Lemma~\ref{lem:two-period-envelope}
(and displayed in the proof of Theorem~\ref{thm:two-period-solution})
are at least $\eta$. Hence $x=y$ is optimal, uniquely when
$\alpha>\alpha_b$.\Halmos
\endproof

\proof{Proof of Theorem \ref{thm:two-period-solution}.}
Set $q=\bar q_2$, $s=\bar x_2$, $\eta=\Phi_2-K_2<0$.
For $K_2>\beta K_1$, differentiation of
\eqref{eq:two-period-envelope} gives slopes
$\Phi_2-K_2\mu/q$ ($x<s,\mu\le q$), $\eta$
($x<s,\mu>q$), and $\Phi_2-K_2\mu/U$ ($x>s$).
Their changes at $\mu=\rho_2q,\rho_2U$ give the first target;
its worst law uses $\{0,q\}$ only when $x^*<s$.
For $K_2\le\beta K_1$, the four slopes, on
$(x<q,\mu\le q)$, $(x<q,\mu>q)$,
$(x>q,\mu\le q)$, and $(x>q,\mu>q)$, respectively, are
\[
\Phi_2-K_2\mu/q,\quad \eta,\quad
\Phi_2-\beta\Phi_1\mu/q,\quad
\eta+(b+h-c)(U-\mu)/(U-q).
\]
The third is positive at $\mu=q$; the fourth changes sign at
$\mu_c=q+[h+\alpha(c-w)](U-q)/(b+h-c)$.
Hence the largest minimizer moves from $0$ to $q$ at $\rho_2q$,
and from $q$ to $U$ at $\mu_c$. Impose $x\ge y$ using
Proposition~\ref{prop:martingale-base-stock}; Lemma~\ref{lem:moment-envelope}
gives the worst-case laws and $\hat V^2=-cy+B_2(x^*,\mu)$.
\Halmos\endproof

\section{Finite-Horizon Characterization}
\label{app:finite-proofs}

Throughout this section, Assumption~\ref{ass:cost-region}, the regular
penalty condition \eqref{eq:regular-penalty}, and $\alpha<\alpha_b$ hold.

\begin{lemma}[Automatic grid ordering]
\label{lem:grid-ordering}
The grids in \eqref{eq:demand-contact-grid}--\eqref{eq:inventory-contact-grid}
satisfy \eqref{eq:automatic-grid} at every finite horizon.
\end{lemma}

\proof{Proof.}
Since $K_t-\Phi_t=b-c-\alpha(c-w)>0$ for $t\ge2$,
$0<\rho_t<1$ and the demand contacts increase strictly.
Regularity $K_2>\beta K_1$ gives
$0<s_{2,1}<q_{2,1}<U$. The threshold recursion preserves
positive differences; the new threshold
$s_{t,t-1}=\beta\Phi_{t-1}U/K_t<U$ because
$K_t-\beta\Phi_{t-1}=b+h-c>0$. Hence positive inventory
thresholds form a strictly ordered suffix. For $t\ge3$,
$K_t>\beta K_{t-1}$. Every inherited positive threshold satisfies
\[
s_{t,j}=\frac{\Phi_{t-1}}{K_t}(s_{t-1,j}-\alpha U)
<\frac{\beta\Phi_{t-1}}{K_t}q_{t-1,j}
<\frac{\Phi_{t-1}}{K_{t-1}}q_{t-1,j}=q_{t,j},
\]
by induction and $q_{t-1,j}\le U$; the newest satisfies
$s_{t,t-1}/q_{t,t-1}=\beta K_{t-1}/K_t<1$.
\Halmos\endproof

\begin{lemma}[Secant recursion]
\label{lem:secant-recursion}
Assume the $(t-1)$-period closed form. For active contacts define
$\Gamma_{t,j}(x)=[F_t^x(q_{t,j})-F_t^x(0)]/q_{t,j}$.
For $a_t\leq j<t$ and $0\leq x\leq q_{t,j}$,
\begin{equation}
\Gamma_{t,j+1}(x)-\Gamma_{t,j}(x)
=K_t\left(\frac1{q_{t,j}}-\frac1{q_{t,j+1}}\right)(x-s_{t,j}).
\label{eq:cleared-secant}
\end{equation}
Moreover, $F_t^x(0)$ is affine on $[0,U]$, and
$F_t^x(q_{t,j})$ is affine on $[0,q_{t,j}]$ for every active $j$.
\end{lemma}

\proof{Proof.}
Write $P=\Phi_{t-1}$, $M=K_{t-1}$ and
$A=(c+(t-1)h-P)U$. At zero demand, the continuation target is zero, so
\begin{equation}
F_t^x(0)=A+(P-c)\alpha U+\{h+\beta(P-c)\}x.
\label{eq:zero-contact-affine}
\end{equation}
For active $j<t$, the inventory recursion gives
$s_{t-1,j}>\alpha U$. At the mean
$q_{t,j}=\rho_{t-1}q_{t-1,j}$, the continuation target is
$s_{t-1,j}$. Thus, when $x\leq q_{t,j}$, the continuation inventory
$\alpha(U-x)$ does not bind that target.
Set
\[
\zeta_j=(t-2)c-(t-1)h
+M\left\{\frac{s_{t-1,j}}{q_{t-1,j}}
+\sum_{k=j+1}^{t-1}\frac{s_{t-1,k}-s_{t-1,k-1}}{q_{t-1,k}}\right\},
\quad j<t,
\]
and $\zeta_t=(t-2)c-(t-1)h+P$.
Substitution of the continuation value gives, for every active $j$,
\begin{equation}
F_t^x(q_{t,j})=bq_{t,j}+A+q_{t,j}\zeta_j
-c\alpha U+(-b+\alpha c)x,
\qquad x\leq q_{t,j}.
\label{eq:positive-contact-affine}
\end{equation}
Since $K_t=b+h-c+\beta P$, subtraction gives
\(\Gamma_{t,j}(x)=b+\zeta_j-\frac{P\alpha U+K_t x}{q_{t,j}}.\)
For consecutive inherited contacts,
\[
\zeta_{j+1}-\zeta_j
=-P s_{t-1,j}\left(\frac1{q_{t,j}}-\frac1{q_{t,j+1}}\right).
\]
The inventory recursion now yields \eqref{eq:cleared-secant}.
The last pair follows by substituting $\zeta_t$ and
$s_{t,t-1}=\beta PU/K_t$.\Halmos
\endproof

\begin{lemma}[Global adjacent-contact dominance]
\label{lem:automatic-contact-dominance}
Suppose Theorem~\ref{thm:finite-horizon} holds at horizon $t-1$.
For $j=1,\ldots,t$, extend the secants to all demand contacts by
\[
\Gamma_{t,j}(x):=\frac{F_t^x(q_{t,j})-F_t^x(0)}{q_{t,j}}.
\]
Then, for every $x\in[0,U]$,
\begin{equation}
\Gamma_{t,j+1}(x)\geq\Gamma_{t,j}(x)
\qquad\text{whenever }x\geq\max\{0,s_{t,j}\}.
\label{eq:global-contact-dominance}
\end{equation}
In particular, all demand contacts discarded by $J_t(x)$ lie below
the chord joining $0$ and $q_{t,J_t(x)}$.
\end{lemma}

\proof{Proof.}
Put $n=t-1$, $P=\Phi_n$, $p_j=q_{n,j}$,
$z_j=(s_{n,j})_+$ ($z_0=0,z_n=U$),
$m_j=\rho_np_j=q_{t,j}$ ($j\le n$), $m_{n+1}=U$,
and $\vartheta_j=m_j/m_{j+1}$. The inductive value implies
\begin{equation}
\hat V^n(y,m_j)=C_j-cy+PH_j(y),\qquad
\hat V^n(y,0)=C_0+(P-c)y,
\label{eq:contact-value-decomposition}
\end{equation}
where $C_j=\hat V^n(0,m_j)$, $C_0=\hat V^n(0,0)$, and
\[
H_j(y)=\sum_{k=j+1}^{n}\left(1-\frac{p_j}{p_k}\right)
\{(y-z_{k-1})_+-(y-z_k)_+\},\quad H_{n+1}=0.
\]
Each $H_j$ is nondecreasing, $1$-Lipschitz, and zero for $y\le z_j$.
The value formula gives
\begin{align}
H_j(y)&=(1-\vartheta_j)(y-z_j)_+
  +\vartheta_jH_{j+1}(y),\label{eq:contact-H-recursion}\\
C_j&=(1-\vartheta_j)C_0+\vartheta_jC_{j+1}
  +P(1-\vartheta_j)z_j.\label{eq:contact-C-recursion}
\end{align}
We also have
\begin{equation}
z_j\le\alpha U+\beta m_{j+1} \quad(j\le n).
\label{eq:contact-lookahead}
\end{equation}
Equality holds for $j=n$. For $j<n$, the threshold recursion gives
$z_j\le\beta\lambda_nq_{n,j+1}$, while
\[
\beta(\lambda_n-\rho_n)_+
=\frac{\beta^n(\alpha w-h)_+}{K_n}\le\alpha.
\]
If $\alpha w>h$, the final inequality follows from
$K_n=(K_n-\Phi_n)+c-w+h/\alpha
+\beta^n(w-h/\alpha)>\beta^n(w-h/\alpha)$;
otherwise it is immediate. Since $q_{n,j+1}\le U$, this proves
\eqref{eq:contact-lookahead}.

Fix $j\le n$ and write $a=m_j$, $d=m_{j+1}$,
$\vartheta=a/d$, $z=z_j$, $s=\alpha(U-x)$. Set
\[
E_j(x)=(1-\vartheta)F_t^x(0)+\vartheta F_t^x(d)-F_t^x(a)
=a\{\Gamma_{t,j+1}(x)-\Gamma_{t,j}(x)\}.
\]
Using \eqref{eq:contact-value-decomposition}--
\eqref{eq:contact-C-recursion}, for $x\le a$,
\[
E_j(x)=(1-\vartheta)
\bigl[(b+h-c)x+P\{s+x-z-(s-z)_+\}\bigr].
\]
If $s\ge z$ it is nonnegative; otherwise it equals
$(1-\vartheta)K_t(x-r_j)$, where
$r_j=P(z-\alpha U)/K_t=s_{t,j}<a$ whenever $r_j>0$.
For $a\le x\le d$, put $y=s+x-a$, $\delta=\vartheta(d-x)$:
\[
E_j(x)=(b+h-c)\delta+PQ,\quad
Q=(1-\vartheta)\{s+x-z-(y-z)_+\}
+\vartheta\{H_{j+1}(s)-H_{j+1}(y)\}.
\]
For $y\ge z$, Lipschitz continuity gives $Q\ge\delta$.
For $y<z$, $E_j$ is affine and nonnegative at the bounding
points $x=a$, $y=z$, and $x=d$, where
$E_j(d)=P(1-\vartheta)(\beta d+\alpha U-z)\ge0$.
Finally, for $x\ge d$ write $y_a=s+x-a$, $y_d=s+x-d$:
\[
E_j(x)=P\bigl[(1-\vartheta)\{s+x-z-(y_a-z)_+\}
+\vartheta\{H_{j+1}(y_d)-H_{j+1}(y_a)\}\bigr].
\]
This is nonnegative by the $1$-Lipschitz bound when $y_a\ge z$;
otherwise it is $P(1-\vartheta)(\beta x+\alpha U-z)\ge0$.
Thus $E_j(x)\ge0$ for $x\ge\max\{0,s_{t,j}\}$;
telescoping proves the discarded-contact claim.\Halmos
\endproof

\begin{lemma}[Automatic finite-horizon worst-case expectation]
\label{lem:finite-envelope-validity}
Under the $(t-1)$-period value representation of
Theorem~\ref{thm:finite-horizon}, the least concave majorant of
$F_t^x$ is its polygonal interpolant at
$0,q_{t,J_t(x)},\ldots,q_{t,t}=U$, for every $x\in[0,U]$.
Its value at $\mu$ is the corresponding worst-case expectation.
\end{lemma}

\proof{Proof.}
Fix $x$ and $J=J_t(x)$. Lemma~\ref{lem:automatic-contact-dominance}
puts discarded contacts below the origin chord ending at $q_{t,J}$.
For $J<t$, $x<s_{t,J}<q_{t,J}$, so
Lemma~\ref{lem:secant-recursion} applies to retained contacts.
For successive chord slopes $m_{t,j}$ it gives
\[
 m_{t,J}-\Gamma_{t,J}=\frac{K_t}{q_{t,J}}(x-s_{t,J})\le0,
 \qquad m_{t,j+1}-m_{t,j}
 =-\frac{K_t}{q_{t,j+1}}(s_{t,j+1}-s_{t,j})<0.
\]
The polygon is concave; for $J=t$ only the endpoint chord remains.
Between consecutive full-grid contacts, partition further at $D=x$
and crossings of $R(x,D)$ with continuation thresholds. On each
smooth segment $F_t^x$ has second derivative $0$ or
$2K_{t-1}/q_{t-1,j}\ge0$; inventory-threshold switches have
nonnegative derivative jumps by the induction slope ordering, and
the jump at $D=x$ is at least $b+h-c>0$ by
\eqref{eq:inventory-slope}.
Thus each full-grid chord dominates its interval. The polygon
therefore majorizes $F_t^x$ globally and touches every retained
contact, proving it is the least concave majorant.\Halmos
\endproof

\proof{Proof of Theorem \ref{thm:finite-horizon}.}
We use induction on $t$, starting from Proposition~\ref{prop:terminal-solution}.
Suppose the stated value representation holds at horizon $t-1$.
Lemmas~\ref{lem:grid-ordering} and~\ref{lem:finite-envelope-validity}
give the exact worst-case expectation of $F_t^x$ for every $x\in[0,U]$;
in particular, the result does not require an additional contact condition.

For $x<U$, let $J=J_t(x)$. On the interior of each inventory interval,
the slope of the fixed-action objective is
\[
\frac{\partial B_t(x,\mu)}{\partial x}
=\Phi_t-K_t\min\{\mu/q_{t,J},1\}.
\]
For the origin chord this follows from
\eqref{eq:zero-contact-affine}--\eqref{eq:positive-contact-affine};
for a chord between positive contacts the slope is
$\Phi_t-K_t=c+\alpha(c-w)-b<0$.
Since $q_{t,J}$ increases with the inventory interval, these slopes are
nondecreasing in $x$. Their zero crossings occur at
$\mu=\rho_tq_{t,J}$. The largest unconstrained minimizer is therefore
$\chi_t(\mu)$ in \eqref{eq:finite-target}, and imposing $x\geq y$
gives $x_t^*(y,\mu)=\max\{y,\chi_t(\mu)\}$.

At capacity, the envelope is the endpoint chord and
\(B_t(U,\mu)=(c+th)U+((t-1)c-th)\mu.\)
For an optimal $x_t^*<U$, its inventory interval satisfies
$\mu<\rho_tq_{t,J_t(x_t^*)}$, so the relevant chord is the origin
chord. The same holds on every inventory interval to its right.
Integrating the action slopes from $x_t^*$ to $U$ gives
\eqref{eq:finite-fixed-action-value} on $[\chi_t(\mu),U]$, hence
the stated value formula.
This closes the induction.

For $x_t^*<U$, the target formula also implies
$\mu<\rho_tq_{t,J_t(x_t^*)}<q_{t,J_t(x_t^*)}$.
The grid ordering gives $x_t^*<q_{t,J_t(x_t^*)}$;
when $x_t^*=U$, the selected contact is $q_t^*=U\geq\mu$.
Thus $q_t^*\geq\max\{x_t^*,\mu\}$ and the origin chord attains
the envelope at the optimal action. Lemma~\ref{lem:moment-envelope}
gives the stated two-point distribution. The target and contact
selections are Borel, and Proposition~\ref{prop:dp-equivalence}
composes the conditional selections into an attaining joint law.

Every positive realization under that law exhausts ready inventory,
leaving $R(x_t^*,q_t^*)=\alpha(U-x_t^*)$.
A zero demand realization is absorbing since
$\mathcal M(0)=\{\delta_0\}$.
Before the first zero, each state has a unique continuing positive
branch under the selected deterministic policy and conditional laws.
Hence the first zero's date, or its absence, determines the complete
path, giving at most $T+1$ demand paths.\Halmos
\endproof

\proof{Proof of Proposition \ref{prop:threshold-stabilization}.}
Set $\lambda_t=\Phi_{t-1}/K_t
=\Phi_{t-1}/(b+h-c+\beta\Phi_{t-1})$.
It increases in $\Phi_{t-1}$ and tends to
$\lambda_\infty\in(0,1)$: $\lambda_t$ increases for
$\alpha<\alpha_h$, decreases for $\alpha>\alpha_h$,
and is constant at equality. For age $\ell$, put
$r_{t,\ell}=s_{t,t-\ell}/U$. The recursions
$r_{t,1}=\beta\lambda_t$ and
$r_{t,\ell}=\lambda_t(r_{t-1,\ell-1}-\alpha)$ give
\begin{equation}
r_{t,\ell}>0\quad\Longleftrightarrow\quad
\sum_{k=1}^{\ell-1}
\left(\prod_{m=1}^k\lambda_{t-\ell+m}\right)^{-1}
<\frac\beta\alpha.
\label{eq:threshold-age-test}
\end{equation}
Positive ages form an initial block. Their limit test is
$\lambda_\infty^{\ell-1}>\alpha/(1-\beta\lambda_\infty)$,
equivalently $\ell-1<\zeta$. For $\alpha\le\alpha_h$,
the finite sums decrease to the limit from above; the count rises to
$N_\infty=\lceil\zeta\rceil$. For $\alpha>\alpha_h$, the
sums increase to the limit from below, so an integer-boundary
age $\zeta$ remains positive at each finite horizon but vanishes
in the limit. Eventually $N_\infty=1+\lfloor\zeta\rfloor$:
all younger ages persist and age $N_\infty+1$ disappears.
\Halmos\endproof

\section{Recovery-Dominant Regimes III and IV}
\label{app:recovery-dominant}

\proof{Proof of Theorem \ref{thm:recovery-dominant-solution}.}
Put
\[
\eta=c+\alpha(c-w)-b=(c-w)(\alpha-\alpha_b)\ge0,\quad
\nu_t=c-\eta\frac{1-(-\alpha)^{t-1}}{1+\alpha},\quad
M_t=\Phi_t-c+\nu_t.
\]
Then $\nu_1=c$, $\nu_t=b+\alpha w-\alpha\nu_{t-1}$,
$0<\nu_t\le c$ and $M_t>0$: if $h\ge\alpha w$ then
$M_t\ge\nu_t$; otherwise $\Phi_t-c>(h-\alpha w)/\alpha$
and $\nu_t\ge c-\eta>w$.
Let $q=\rho_1U$ and
\[
\begin{gathered}
v_{t,0}(y)=thU-(\Phi_t-c)(U-y),\quad
v_{t,U}(y)=U\sum_{r=1}^t\nu_r-\nu_ty,\\
s_t^\circ=\frac U{M_t}
\left(\Phi_1-\alpha\sum_{r=1}^{t-1}M_r\right)_{+},\qquad
\tau(\mu)=\min\left\{\frac{(1-\rho_1)\mu}{q},1-\frac\mu U\right\}.
\end{gathered}
\]
Regularity gives $s_2^\circ=\bar x_2<q$, and the recursion
$s_t^\circ=M_{t-1}(s_{t-1}^\circ-\alpha U)_{+}/M_t$
keeps $0\le s_t^\circ<q$. Induction gives
\begin{equation}
\hat V^t(y,\mu)
=\left(1-\frac\mu U\right)v_{t,0}(y)
 +\frac\mu Uv_{t,U}(y)
 +M_t(s_t^\circ-y)_{+}\tau(\mu).
\label{eq:recovery-dominant-value}
\end{equation}
The base is Proposition~\ref{prop:terminal-solution}.
Inductively $F_t^x$ is convex on $[0,q]$, $[q,U]$; for $x\le q$,
the excess at $q$ above the endpoint chord equals
\[
M_{t-1}(s_{t-1}^\circ-\alpha U+\alpha x)_{+}
-(M_t+\alpha M_{t-1})x,
\]
whose positive part is $M_t(s_t^\circ-x)_{+}$; for $x\ge q$
that excess is convex and nonpositive at its endpoints.
Consequently, the least concave majorant uses adjacent chords
through $0,q,U$ only for $x<s_t^\circ$; otherwise the endpoint
chord dominates. Its three possible action slopes $\Phi_t-M_t\mu/q$, $\Phi_t-M_t$ and
$\Phi_t-M_t\mu/U$ are at least $c-\nu_t\ge0$; hence
$x=y$ is optimal nonterminally (unique for $\alpha>\alpha_b$).
Substitution proves \eqref{eq:recovery-dominant-value}, with
adjacent-contact laws attaining it. Since $\alpha>0$ and $M_t$
is bounded away from zero, $s_t^\circ$ reaches zero in finite time;
the terminal rule follows from Proposition~\ref{prop:terminal-solution}.\Halmos
\endproof

\section{Approximation Policies for the Complementary Case}
\label{app:complementary-approx}
\label{subsec:complementary-approx}

Under Assumption~\ref{ass:cost-region}, the complementary condition
$K_2\le\beta K_1$ occurs only in Regimes II/IV, with
$h<w/4$ and $0<b-c\le w-2\sqrt{hw}$. The two-period contact
$\bar q_2$ need not persist, so we bound implementable policies using
\eqref{eq:finite-coefficients}.

In Regime II, define the cost-based threshold
\begin{equation*}
\theta:=\frac{h+\alpha(c-w)}{b+h-c}\in(0,1),
\qquad \underline\rho:=\min\{\theta,\rho_1\}.
\end{equation*}
Use the policy
\begin{equation}
\tilde x_t(y,\mu)=
\begin{cases}
x_t^*(y,\mu),&t=1,2,\\[1mm]
y,&t\geq3,\ \mu<\theta U,\\[1mm]
U,&t\geq3,\ \mu\geq\theta U.
\end{cases}
\label{eq:regime-II-approx-policy}
\end{equation}
Relative to capacity continuation, action $x$ changes cost by
$(U-x)[(b+h-c)\mu/U-h-\alpha(c-w)]$; hence $\theta U$
is the benchmark (not necessarily optimal) switching mean.

In Regime IV, use no replenishment in every nonterminal period:
\begin{equation}
\tilde x_t(y,\mu)=
\begin{cases}
x_1^*(y,\mu),&t=1,\\[1mm]
y,&t\geq2.
\end{cases}
\label{eq:regime-IV-approx-policy}
\end{equation}
Let $\tilde V^t$ denote the policy's worst-case cost and
$g:=\beta K_1-K_2\ge0$.

\begin{proposition}[Finite-horizon approximation guarantees]
\label{prop:complementary-approx}
Under Assumption~\ref{ass:cost-region}, let $K_2\leq\beta K_1$.
Both policies are exact at $t=1,2$ and $\mu\in\{0,U\}$.
For any state $(y,\mu)\in[0,U]^2$:
\begin{enumerate}
\item In Regime II, define the dimensionless constant
\begin{equation}
\mathcal K_{\mathrm{II}}:=\frac{\beta\Phi_2(1-\underline\rho)}
 {h(1-\underline\rho)+c\underline\rho}.
\label{eq:approx-II-kappa}
\end{equation}
For every $t\geq3$,
\begin{equation}
1\leq\frac{\tilde V^t(y,\mu)}{\hat V^t(y,\mu)}
\leq1+\frac{\mathcal K_{\mathrm{II}}}{t},
\qquad
0\leq\tilde V^t(y,\mu)-\hat V^t(y,\mu)
\leq\beta\Phi_2(1-\underline\rho)U.
\label{eq:regime-II-regret-bound}
\end{equation}
\item In Regime IV, for every $t\geq1$,
\begin{equation}
0\leq\tilde V^t(y,\mu)-\hat V^t(y,\mu)
\leq g\min\{(1-\rho_1)\mu,\rho_1(U-\mu)\}
\leq\frac{gU}{4}.
\label{eq:regime-IV-regret-bound}
\end{equation}
In particular, the policy is optimal at every horizon if $g=0$.
\end{enumerate}
\end{proposition}

Both policies avoid the demand grid. The Regime II relative error
is $O(1/t)$; the Regime IV gap is bounded uniformly in $t$.

For Regime IV's relative bound, define
\begin{equation}
\eta:=c+\alpha(c-w)-b\geq0,\qquad
\nu_r:=c-\eta\frac{1-(-\alpha)^{r-1}}{1+\alpha}>0.
\label{eq:approx-IV-boundary-coefficients}
\end{equation}
For $t\geq2$,
\begin{equation}
1\leq\frac{\tilde V^t(y,\mu)}{\hat V^t(y,\mu)}
\leq1+\frac{g\rho_1(1-\rho_1)}
 {th(1-\rho_1)+\rho_1\sum_{r=1}^{t-1}\nu_r}.
\label{eq:approx-IV-relative-bound}
\end{equation}
The excess is $O(1/t)$.

\proof{Proof of Proposition~\ref{prop:complementary-approx}(1).}
Put $p=\mu/U$, $a=h+\alpha(c-w)$,
$d=b-c-\alpha(c-w)>0$, $L=b+h-c=a+d$, and
\begin{equation}
\Lambda(\mu)=h(U-\mu)+c\mu,\quad
C_t^{\mathrm{cap}}(y,\mu)=t\Lambda(\mu)+c(U-y)-c\mu.
\label{eq:approx-capacity-cost}
\end{equation}
Capacity replenishment incurs this cost. The absorbing endpoint law
$\{0,U\}$, together with
$\hat V^t(y,0)=thU-(\Phi_t-c)(U-y)$ and
$\hat V^t(y,U)=tcU-cy$, gives first-action slope
$\Phi_t-K_tp$. Hence
\begin{equation}
\hat V^t(y,\mu)\ge C_t^{\mathrm{cap}}(y,\mu)
 -(U-y)(\Phi_t-K_tp)_{+}.
\label{eq:approx-endpoint-improvement}
\end{equation}
For $W_0(y,\mu)=c(U-y)-c\mu$, the stage identity is
\begin{equation}
\mathbb E[C(y,x,D)+W_0(R(x,D),D)]
 -\Lambda(\mu)-W_0(y,\mu)
 =L\mathbb E(D-x)_{+}-a(U-x).
\label{eq:approx-BR-stage-identity}
\end{equation}
The endpoint law maximizes the right side at $(U-x)(Lp-a)$.
Stagewise comparison with capacity, using exact terminal actions, gives
$\tilde V^t\le C_t^{\mathrm{cap}}-(U-y)(a-Lp)_+$.
For a second comparison set
\begin{equation}
\tilde\lambda_1=\max\{K_1,\Phi_1/\theta\},\qquad
\tilde\lambda_t=\frac L\alpha+
 (\tilde\lambda_1-L/\alpha)\beta^{t-1}.
\label{eq:approx-BR-comparison-coefficients}
\end{equation}
Then $\tilde\lambda_t=L+\beta\tilde\lambda_{t-1}$,
$\Phi_t=a+\beta\Phi_{t-1}$ and
$r_t=\Phi_t/\tilde\lambda_t\le\theta$. We show inductively
\begin{equation}
\tilde V^t\le C_t^{\mathrm{cap}}-(U-y)
\max\{(a-Lp)_+,(\Phi_t-\tilde\lambda_tp)_+\},\qquad t\ge2.
\label{eq:approx-II-sandwich}
\end{equation}
At $t=1$, the exact terminal value gives the comparison with
$\tilde\lambda_1\ge K_1$; at $t=2$ the policy is optimal and is no
worse than either the threshold or capacity comparison. For $t\ge3$
and $p<\theta$, the policy chooses $x=y$. The inductive continuation
bound makes its excess at most
$\sup_{Q\in\mathcal M(\mu)}\mathbb E_Qq(D)$, where
\[
q(D)=L(D-y)_+-a(U-y)
 -(U-R(y,D))(\Phi_{t-1}-\tilde\lambda_{t-1}D/U)_+.
\]
For $r=r_{t-1}$, $q$ is convex on $[0,rU]$, $[rU,U]$
(with nonnegative quadratic coefficients and derivative jumps), and
$q(0)=-(U-y)\Phi_t$, $q(rU)\le(U-y)(Lr-a)$,
$q(U)=(U-y)d$. The line joining the bounds at $0,rU$ is
$(U-y)(-\Phi_t+\tilde\lambda_tD/U)$; the line through $rU,U$
is $(U-y)(-a+LD/U)$. They agree at $D=rU$, since
$\Phi_t=a+\beta\Phi_{t-1}$ and
$\tilde\lambda_t=L+\beta\tilde\lambda_{t-1}$.
Their slopes decrease ($\tilde\lambda_t>L$). Convexity on the two
intervals and the three contact bounds therefore yield the
\emph{global} majorant
\begin{equation}
q(D)\le(U-y)\min\{-\Phi_t+\tilde\lambda_tD/U,
 -a+LD/U\}.
\label{eq:approx-BR-demand-majorant}
\end{equation}
Because the minimum of these two affine functions is concave,
Jensen's inequality bounds the worst-case expectation by its value
at $\mu$. For $p\ge\theta$ the policy selects capacity and both
improvements in \eqref{eq:approx-II-sandwich} vanish.

Put $A_t=\beta\Phi_{t-1}$ and
$\tilde\delta_t=\tilde\lambda_t-K_t
=\beta(\tilde\lambda_{t-1}-\Phi_{t-1})>0$.
Weighted averaging gives $\underline\rho=r_1\le r_{t-1}
\le r_t\le\theta$ and $A_t/(A_t+\tilde\delta_t)=r_{t-1}$;
thus
\begin{equation}
0\le\tilde V^t-\hat V^t\le(U-y)
\min\{\tilde\delta_tp,A_t(1-p),(\Phi_t-K_tp)_+\}.
\label{eq:approx-II-state-bound}
\end{equation}
The first two terms meet at $p=r_{t-1}$, so their minimum is at
most $A_t(1-r_{t-1})\le\beta\Phi_2(1-\underline\rho)$.
For the relative bound only $p<\rho_t$ matters. With $v=U-y$,
\eqref{eq:approx-endpoint-improvement} lower-bounds the value by
$U\{th(1-p)+(t-1)cp\}+v(c-\Phi_t+K_tp)>0$.
The ratio increases in $v$, so take $v=U$ and set
$\mathcal D_t(p)=th+c-\Phi_t+\{(t-1)c-th+K_t\}p$.
Since $\mathcal D_t(0)>0$ and $\mathcal D_t(1)=tc+d>0$,
$\tilde\delta_tp/\mathcal D_t(p)$ increases and
$A_t(1-p)/\mathcal D_t(p)$ decreases. Their minimum is largest
at $r=r_{t-1}$, giving
\[
\frac{\tilde V^t-\hat V^t}{\hat V^t}
\le\frac{A_t(1-r)}{(1-r)(th+c-\Phi_t)+r(tc+d)}
\le\frac{\beta\Phi_2(1-\underline\rho)}
 {t\{h(1-\underline\rho)+c\underline\rho\}}
=\frac{\mathcal K_{\mathrm{II}}}{t}.
\]
Here $\Phi_t<c$, $r\ge\underline\rho$ and
$A_t\le\beta\Phi_2$. Endpoint errors vanish and $t=1,2$
are exact.\Halmos\endproof

\proof{Proof of Proposition~\ref{prop:complementary-approx}(2).}
Use $\eta,\nu_t$ from \eqref{eq:approx-IV-boundary-coefficients},
$M_t=\Phi_t-c+\nu_t$, and $v_{t,0},v_{t,U}$ from the proof of
Theorem~\ref{thm:recovery-dominant-solution}. The inequalities
$0<\nu_t\le c$, $M_t>0$ and
$M_t-\beta M_{t-1}=K_1-\nu_{t-1}>0$ do not use regularity.
Set
\[
V_{\mathrm{bd}}^t(y,\mu)
=(1-\mu/U)v_{t,0}(y)+(\mu/U)v_{t,U}(y),\quad
\tau(\mu)=\min\{(1-\rho_1)\mu/\bar q_2,1-\mu/U\}.
\]
The complementary two-period value is
$\hat V^2=V_{\mathrm{bd}}^2+\Omega_2(y)\tau(\mu)$, with
\[
\Omega_2(y)=\begin{cases}
\beta\Phi_1U-K_2y,&y\le\bar q_2,\\
 g\rho_1(U-y)/(1-\rho_1),&y\ge\bar q_2.
\end{cases}
\]
For $\underline V^2=V_{\mathrm{bd}}^2
+K_2(\bar q_2-y)_+\tau(\mu)$, therefore
\begin{equation}
0\le\hat V^2-\underline V^2\le g\rho_1U\tau(\mu).
\label{eq:approx-IV-terminal-perturbation}
\end{equation}
Extend $\underline V^2$ by Bellman recursion and define
$\underline s_2=\bar q_2$,
$\underline s_t=(M_{t-1}/M_t)(\underline s_{t-1}-\alpha U)_+$;
thus $0\le\underline s_t\le\bar q_2$. Induction gives
\begin{equation}
\underline V^t=V_{\mathrm{bd}}^t
+M_t(\underline s_t-y)_+\tau(\mu),\qquad t\ge2.
\label{eq:approx-IV-lower-value}
\end{equation}
Indeed the induction integrand is convex in demand on each of
$[0,\bar q_2]$, $[\bar q_2,U]$: second derivatives on smooth
pieces and all derivative jumps are nonnegative. For $x\le\bar q_2$
the excess at $\bar q_2$ over the endpoint chord is
\[
M_{t-1}(\underline s_{t-1}-\alpha U+\alpha x)_+
-(M_t+\alpha M_{t-1})x,
\]
whose positive part is $M_t(\underline s_t-x)_+$.
For $x\ge\bar q_2$ the excess is convex and nonpositive at both
ends, so the envelope uses $0,\bar q_2,U$ when
$x<\underline s_t$ and $0,U$ otherwise. Its action slopes
$\Phi_t-M_t\mu/\bar q_2$, $\Phi_t-M_t$,
$\Phi_t-M_t\mu/U$ are at least $c-\nu_t\ge0$;
thus $x=y$ minimizes the auxiliary problem and proves
\eqref{eq:approx-IV-lower-value}.
At horizon two, \eqref{eq:approx-IV-terminal-perturbation}
gives $\underline V^2\le\hat V^2=\tilde V^2
\le\underline V^2+g\rho_1U\tau(\mu)$.
The auxiliary recursion has optimal no-order actions; the actual
no-order policy obeys the same recursion above this terminal base.
The two Bellman operators are monotone and invariant under additive
constants, and concavity of $\tau$ gives
$\mathbb E_Q[\tau(\widetilde D)]\le\tau(\mu)$ for
$Q\in\mathcal M(\mu)$. Induction therefore preserves
\[
\underline V^t\le\hat V^t\le\tilde V^t
\le\underline V^t+g\rho_1U\tau(\mu),
\]
proving \eqref{eq:regime-IV-regret-bound}. Moreover
$\hat V^t\ge V_{\mathrm{bd}}^t(y,\mu)
\ge th(U-\mu)+\mu\sum_{r=1}^{t-1}\nu_r>0$ for $t\ge2$.
The error divided by this affine lower bound is maximized at
$\mu=\bar q_2$, giving \eqref{eq:approx-IV-relative-bound};
its denominator grows linearly since
$\nu_r\to(b+\alpha w)/(1+\alpha)>0$. The endpoint and
$g=0$ claims follow.\Halmos\endproof

\section{Proofs for Directional Demand Information}
\label{app:directional-drift}

For $\diamond\in\{\uparrow,\downarrow\}$, set
$F_{\diamond,t}^x(D)=b(D-x)_++h(x-D)_+
+\hat V_\diamond^{t-1}(R(x,D),D)$, and let
$B_{\diamond,t}(x,\mu)$ add $cx+\alpha w(U-x)$ to its
worst-case directional expectation. Then
$\hat V_\diamond^t=-cy+\min_{x\in[y,U]}B_{\diamond,t}(x,\mu)$.

\proof{Proof of Proposition \ref{prop:directional-dp} and Corollary \ref{cor:directional-monotonicity}.}
Write $H=\operatorname{cav}\varphi$ for a continuous demand integrand.
Partitioning admissible laws by their actual means gives
\[
\sup_{Q\in\mathcal M^\uparrow(\mu)}\mathbb E_Q\varphi
=\max_{m\in[\mu,U]}H(m),\qquad
\sup_{Q\in\mathcal M^\downarrow(\mu)}\mathbb E_Q\varphi
=\max_{m\in[0,\mu]}H(m).
\]
Continuity of $H$ ensures attainment, and
Lemma~\ref{lem:moment-envelope} supplies a law of at most two points
at its chosen mean. These ambiguity correspondences are compact and
continuous (mix with $\delta_U$ or $\delta_0$ for lower
hemicontinuity), so Proposition~\ref{prop:dp-equivalence} applies.
Set inclusion gives the value bounds and monotonicity in $\mu$.\Halmos
\endproof

\proof{Proof of Proposition \ref{prop:directional-base-stock}.}
Inventory convexity gives the constrained largest minimizer
$\max\{y,\chi_t^\diamond(\mu)\}$; measurable selection
follows from the maximum theorem.
\Halmos
\endproof

\proof{Proof of Theorem \ref{thm:submartingale-solution} and
Corollary \ref{cor:submartingale-capacity}.}
For the $\alpha<\alpha_b$ part, set $\theta_t=th-(t-1)c$.
At the absorbing reference state $U$, direct induction gives
\(\hat V_\uparrow^t(y,U)=tcU-cy,\qquad x_{\uparrow,t}^*(y,U)=U.\)
For $t\geq2$, the coefficient of $U-x$ in this deterministic
problem is $b-c-\alpha(c-w)>0$.

We first handle $\theta_t\leq0$. This requires $c>h$ and is
equivalent to $t\geq\lceil c/(c-h)\rceil$.
The feasible all-$U$ demand path gives the lower bound
$tcU-cy+\{b-c-\alpha(c-w)\}(U-x)$.
Under the policy that replenishes to $U$ every period, let
$m_r=\mathbb E[\widetilde D_r]$. Its expected cost is
\(c(U-y)+thU+(c-h)\sum_{r=2}^t m_r-hm_1.\)
Because $\mu\leq m_t\leq\cdots\leq m_1\leq U$, this is at most
$c(U-y)+thU+\{(t-1)c-th\}m_1\leq tcU-cy$.
The bounds coincide, proving the capacity result and the theorem
when $\theta_t\leq0$.

For $\theta_t>0$, we prove the formula by induction. The equivalent
value representation is
\begin{equation}
\hat V_\uparrow^t(y,\mu)
=-cy+C_t^\uparrow
 +\left(\Phi_t-K_t\frac{\mu}{U}\right)_{+}
       (y-\omega_t U)_{+}
 -\theta_t(\mu-\rho_tU)_{+}.
\label{eq:submartingale-induction-value}
\end{equation}
For $t=1$, the terminal loss has an endpoint concave envelope;
maximizing it over means in $[\mu,U]$ gives the stated formula.

For the induction step, write $P=\Phi_{t-1}$,
$M=K_{t-1}$, and $q=P/M$.
Fix $x\geq\omega_t U$, write $\bar x=x/U$, and let
$a=\theta_{t-1}/M-\beta(1-\bar x)$.
Here $\theta_t>0$ and $t\geq2$ imply $h>c/2>w/2$, so
$K_2>\beta K_1$ holds automatically; for $t\geq3$,
$K_t>\beta K_{t-1}$ follows from the coefficient identity in
Section~\ref{sec:finite-horizon}. The inequalities
$K_t>\beta M$, $\theta_{t-1}=\theta_t+c-h$, and $P\geq \theta_{t-1}$
imply $0\leq a\leq q$. Substituting the induction formula at
the demand endpoints gives
\begin{equation}
\begin{aligned}
cx+\alpha w(U-x)+F_{\uparrow,t}^x(U)
&=tcU+(K_t-\Phi_t)(U-x),\\
F_{\uparrow,t}^x(0)-F_{\uparrow,t}^x(U)
&=K_t(x-\omega_t U).
\end{aligned}
\label{eq:submartingale-endpoints}
\end{equation}
For $z\in[0,1]$, the normalized endpoint-chord residual is
\begin{equation*}
\begin{aligned}
E(z)&:=\frac{(1-z)F_{\uparrow,t}^x(0)
                  +zF_{\uparrow,t}^x(U)-F_{\uparrow,t}^x(Uz)}{U}\\
&=(b+h-c)\min\{z(1-\bar x),\bar x(1-z)\}
 -\theta_{t-1}\min\{z(1-q),q(1-z)\}\\
&\quad+Pa(1-z)-M(q-z)_{+}(a-\min\{z,\bar x\})_{+}.
\end{aligned}
\end{equation*}
With $u=1-\bar x$, substitution of $\theta_{t-1}=M(a+\beta u)$ reduces this to
\[
\begin{array}{ll}
z\{M(q-z)+u(K_t-\beta M)\},
 &0\leq z\leq\min\{a,\bar x\},\\
Ma(q-z)+uz(K_t-\beta M),
 &a\leq z\leq\min\{\bar x,q\},\quad a\leq\bar x,\\
u(K_t z-\beta P),&q\leq z\leq\bar x,\\
(1-z)(K_t\bar x-\beta P),&\max\{\bar x,q\}\leq z\leq1.
\end{array}
\]
Each expression is nonnegative because $K_t>\beta M$ and
$K_t\bar x-\beta P\geq b+h-c-\theta_t>0$.
If $\bar x<q$, the remaining interval $[\bar x,q]$ has an
affine residual with nonnegative endpoints. Thus the endpoint
chord majorizes the integrand everywhere.
By the endpoint difference in \eqref{eq:submartingale-endpoints},
that chord is nonincreasing in its mean for $x\geq\omega_t U$.
Its maximum over $[\mu,U]$ is attained at $\mu$,
and hence
\begin{equation}
B_{\uparrow,t}(x,\mu)
=tcU+(K_t-\Phi_t)(U-x)
 +K_t\left(1-\frac{\mu}{U}\right)(x-\omega_t U),
\quad x\geq\omega_t U.
\label{eq:submartingale-fixed-action}
\end{equation}
For $x<\omega_t U$, the feasible law $\delta_U$ instead gives
\[
B_{\uparrow,t}(x,\mu)
\geq tcU+(K_t-\Phi_t)(U-x)
>B_{\uparrow,t}(\omega_t U,\mu).
\]
Therefore minimization only requires
\eqref{eq:submartingale-fixed-action}, whose slope is
$\Phi_t-K_t\mu/U$. Its sign gives the stated target and
\eqref{eq:submartingale-induction-value}, completing the induction.
When $\theta_t>0$, capacity is strictly suboptimal at $(y,\mu)=(0,0)$.
This proves that the capacity horizon is exact under the standing
cost ordering $c>h$.
\Halmos
\endproof

\proof{Proof of the recovery-dominant assertion in Theorem \ref{thm:submartingale-solution}.}
Write $k=b/(b+h)$ and
$\eta=c+\alpha(c-w)-b=(c-w)(\alpha-\alpha_b)\geq0$.
Since $R(0,D)=\alpha U$,
$F_{\uparrow,2}^0(D)=bD+\hat V_\uparrow^1(\alpha U,D)$.
Below $D=\rho_1U$, its slope is $b$ when $\alpha\leq k$
and $2b-(b+h)\alpha>b-h$ when $\alpha>k$;
above $\rho_1U$, the slope is $b-h>0$. Thus this function
is strictly increasing, and $\delta_U$ uniquely maximizes its
expectation under upward ambiguity.

At terminal reference $U$,
$\hat V_\uparrow^1(R(x,U),U)=c\{U-\alpha(U-x)\}$.
Danskin's theorem therefore gives
$\partial_x^+B_{\uparrow,2}(0,\mu)=
 c-\alpha w-b+\alpha c=\eta\ge0$.
Convexity of $B_{\uparrow,2}(\cdot,\mu)$ yields the optimal
selection $x_{\uparrow,2}^*(y,\mu)=y$;
when $\alpha>\alpha_b$, the positive derivative makes it unique.
\Halmos
\endproof

\proof{Proof of Proposition \ref{prop:supermartingale-recovery-suppression}.}
Put $k=\alpha_d=b/(b+h)$,
$z=z_2^\downarrow=U(1-k/\alpha)>0$, and
$\eta=c+\alpha(c-w)-b>0$. For $0\le x\le z$,
$R(x,D)\ge\alpha(U-x)\ge\alpha(U-z)=kU$.
The terminal controller thus never replenishes, and the terminal
downward envelope selects mean zero, since its endpoint slope
$b-(b+h)R(x,D)/U\le0$. Accordingly,
$\hat V_\downarrow^1(R(x,D),D)
=\alpha w\{U-R(x,D)\}+hR(x,D)$.
The function $F_{\downarrow,2}^x$ is convex in demand, with
slopes $\alpha w-2h$ for $D<x$ and $b$ for $D>x$;
their jump is $b+2h-\alpha w>0$.
Its endpoint difference is
$\Delta(x)=bU-(b+2h-\alpha w)x$, positive on $[0,z]$:
indeed, writing $q=\alpha-k\in(0,1-k)$,
\[
\frac{\alpha\Delta(z)}{U}=kb-q(2h-\alpha w),\quad
q(2h-\alpha w)\le h(\alpha-k)(2-\alpha)
\le\frac{h^2}{b+h}<\frac{b^2}{b+h}=kb.
\]
The second inequality uses monotonicity of
$(\alpha-k)(2-\alpha)$ on $[k,1]$.
Hence the worst law on $[0,z]$ is
$(1-\mu/U)\delta_0+(\mu/U)\delta_U$, and the action slope is
\[
\partial_x B_{\downarrow,2}(x,\mu)
=\gamma_0-(\mu/U)(\gamma_0-\gamma_1),\quad
\begin{cases}
\gamma_0=c+(2-\alpha)(h-\alpha w),\\
\gamma_1=c-b-\alpha h-\alpha\beta w.
\end{cases}
\]
Since $\alpha(2-\alpha)\le1$,
$\gamma_0>c-w>0$ and $\gamma_1<c-b<0$;
the slope changes sign at $\mu_2^\dagger$.
At $x=z$, the endpoint chord is strictly above the demand
integrand in $(0,U)$, hence its endpoint law is the unique
maximizer when $0<\mu<U$. To the right, the $D=0$ action slope
is still $\gamma_0$, while $R(x,U)<kU$ triggers terminal
replenishment, giving the $D=U$ slope
$c-\alpha w-b+\alpha c=\eta>0$.
The one-sided Danskin formula (also immediate at the endpoints) gives
\[
\partial_x^+ B_{\downarrow,2}(z,\mu)
=(1-\mu/U)\gamma_0+(\mu/U)\eta>0.
\]
Convexity now yields largest minimizer $0$ for
$\mu<\mu_2^\dagger$ and $z$ for $\mu\ge\mu_2^\dagger$
(the objective is flat on $[0,z]$ at equality). Imposing
$x\ge y$ gives the asserted rule.\Halmos\endproof

\proof{Proof of Theorem \ref{thm:supermartingale-solution}.}
At $t=1$, the endpoint envelope for action $x$ is
$hx+\{b-(b+h)x/U\}m$.
Maximizing over $0\leq m\leq\mu$ selects $m=\mu$ for
$x<bU/(b+h)$ and $m=0$ for $x>bU/(b+h)$. The corresponding
action slopes are $\Phi_1-K_1\mu/U$ and $\Phi_1$, which give
$s_{1,1}^\downarrow=bU/(b+h)$, the stated terminal target, and the
terminal value.

For the induction step, write $n=t-1$, $P=\Phi_n$, and
$z_j=s_{n,j}^\downarrow$. Define the untruncated inherited switches
\(r_j=\frac{P}{K_t}(z_j-\alpha U),\qquad j=1,\ldots,n.\)
Along with the value formula, we prove that the capacity mean slope
$A_r:=(\hat V_\downarrow^{r-1}(0,U)/U-rh)_{+}$ is zero whenever
$s_{r,r}^\downarrow<U$. This holds at $r=1$.
We establish the envelope before minimizing it.

\emph{Ordering.}
The ratios $\rho_s$ are nondecreasing: $\Phi_s$ is nondecreasing,
$K_s-\Phi_s=b-c-\alpha(c-w)>0$ for $s\geq2$, and
$K_1-\Phi_1$ exceeds this difference by $\alpha c$.
Induction gives
\[
0<q_{s,1}<\cdots<q_{s,s}=U,\quad
\frac{q_{s,j}}{q_{s,j+1}}=\rho_j,\quad
0\leq s_{s,1}^\downarrow\leq\cdots\leq s_{s,s}^\downarrow\leq U,
\quad s_{s,j}^\downarrow\leq q_{s,j}.
\]
For inherited switches, $r_j\leq q_{t,j}$ follows from
$K_t\geq K_{t-1}$ for $t\geq3$. At $t=2$, it follows from
\[
K_2-K_1\left(\frac{b}{b+h}-\alpha\right)
=2h-\alpha w+\alpha(b-c)+\alpha^2w>0.
\]
Positive targets are strictly ordered: the recursion preserves positive
inherited differences, and the new highest target strictly exceeds
the last inherited target, which is below $U$.
Monotonicity of $\rho_s$ also implies
\[
\beta q_{t,j+1}+\alpha U\geq z_j,
\qquad 1\leq j\leq n.
\]
The highest target satisfies $s_{s,s}^\downarrow\geq bU/(b+h)$.
Indeed, with $k=b/(b+h)$, its untruncated recursion at $t$ exceeds
$kU$ by at least $U\{kc-\alpha(1-k)P\}/K_t\geq0$;
here $\alpha P\leq h+\alpha(c-w)<bc/h$.
Consequently $s_{s,s}^\downarrow\geq\alpha U$ under the theorem's
additional assumption.

\emph{Adjacent-contact comparisons.}
As in the proof of Lemma~\ref{lem:automatic-contact-dominance},
put $p_j=q_{n,j}$, $m_j=\rho_np_j=q_{t,j}$,
$m_{n+1}=U$, and $\vartheta_j=m_j/m_{j+1}$. The induction value has
\[
\hat V_\downarrow^n(y,m_j)=C_j-cy+PH_j(y),\qquad
\hat V_\downarrow^n(y,0)=C_0+(P-c)y,
\]
where $C_j=\hat V_\downarrow^n(0,m_j)$,
$C_0=\hat V_\downarrow^n(0,0)$, and
\[
H_j(y)=\sum_{\ell=j+1}^{n}
\left(1-\frac{p_j}{p_\ell}\right)
\bigl\{(y-z_{\ell-1})_{+}-(y-z_\ell)_{+}\bigr\}
+(y-z_n)_{+},\qquad H_{n+1}(y)=(y-z_n)_{+}.
\]
The functions $H_j$ are nondecreasing and $1$-Lipschitz, vanish
for $y\le z_j$, and satisfy the same recursions
\eqref{eq:contact-H-recursion}--\eqref{eq:contact-C-recursion}.
For adjacent contacts $a=m_j$, $d=m_{j+1}$, put
$\vartheta=a/d$, $z=z_j$, $s=\alpha(U-x)$, and
\[
E_j(x)=(1-\vartheta)F_{\downarrow,t}^x(0)
+\vartheta F_{\downarrow,t}^x(d)-F_{\downarrow,t}^x(a).
\]
The three-range calculation in the proof of
Lemma~\ref{lem:automatic-contact-dominance} applies unchanged:
for $x\le a$,
\[
E_j(x)=(1-\vartheta)
\bigl[(b+h-c)x+P\{s+x-z-(s-z)_{+}\}\bigr],
\]
which is negative exactly when $r_j>0$ and $x<r_j$.
For $a\le x\le d$, the same Lipschitz bound and
$\beta d+\alpha U\ge z$ show $E_j(x)\ge0$;
for $x\ge d$, the bound follows from the $1$-Lipschitz property
if $s+x-a\ge z$, and from
$E_j(x)=P(1-\vartheta)(\beta x+\alpha U-z)\ge0$
otherwise. Therefore
\begin{equation}
E_j(x)<0\quad\Longleftrightarrow\quad r_j>0\text{ and }x<r_j;
\qquad E_j(x)\ge0\text{ otherwise}.
\label{eq:downward-contact-residual}
\end{equation}

\emph{Concave envelope.}
Let $J$ be the first index with $x<r_J$, or $J=t$ if none exists,
and define
$\Gamma_j(x)=\{F_{\downarrow,t}^x(q_{t,j})
                  -F_{\downarrow,t}^x(0)\}/q_{t,j}$.
Since $E_j(x)=q_{t,j}\{\Gamma_{j+1}(x)-\Gamma_j(x)\}$,
\eqref{eq:downward-contact-residual} shows that every contact below
$q_{t,J}$ lies below its origin chord.
When $J<t$, for the retained contacts $j\geq J$, current demand clears
inventory and continuation inventory lies below the corresponding
target. If $\ell_j(x)$ is the chord slope from $q_{t,j}$ to
$q_{t,j+1}$, the same substitution gives
\[
\ell_J(x)-\Gamma_J(x)=\frac{K_t}{q_{t,J}}(x-r_J)<0,
\qquad
\ell_{j+1}(x)-\ell_j(x)
=-\frac{K_t}{q_{t,j+1}}(r_{j+1}-r_j)\leq0.
\]
Thus the polygon through $0,q_{t,J},\ldots,q_{t,t}$ is concave
(a single chord if $J=t$).
Between consecutive full-grid contacts, partition demand further at
$D=x$ and $R(x,D)=z_j$. On each smooth piece the integrand has
second derivative either zero or $2K_n/q_{n,j}>0$.
Inventory switches produce nonnegative derivative jumps, and the
jump at $D=x$ is at least $b+h-c>0$ by
\eqref{eq:inventory-slope}.
The integrand is therefore convex between full-grid contacts.
The polygon dominates those contacts and hence the entire integrand;
because it touches every retained contact, it is the least concave
majorant.

\emph{Minimization and value.}
On a retained positive-contact face, each contact clears inventory
and has continuation inventory below its target. The fixed-action
slope is therefore $c+\alpha(c-w)-b<0$.
An action for which the downward envelope selects such a face
cannot be optimal. On an increasing origin face ending at $q_{t,j}$,
the slope is $\Phi_t-K_t\mu/q_{t,j}$; once the prefix maximum
selects mean zero, the slope is $\Phi_t>0$.

To locate the last change, first suppose $z_n<U$. The induction
capacity anchor has zero mean slope in this case, and endpoint
substitution on the last inventory face gives
\[
F_{\downarrow,t}^x(U)-F_{\downarrow,t}^x(0)
=K_t\left\{\frac{P}{K_t}(z_n-\alpha U)
                       +\frac{b}{K_t}U-x\right\}.
\]
This yields the stated highest target, clipped at capacity.
If $z_n=U$, the same untruncated candidate is at least $U$, since
$(\beta P+b-K_t)U/K_t=(c-h)U/K_t\geq0$; the highest target
therefore remains at capacity.
The inherited changes occur at the $r_j$, clipped at zero and $U$.
Their ordering and the signs of the origin-face slopes give the
stated switching means $\rho_t q_{t,j}$ and the largest minimizer
$\chi_t^\downarrow(\mu)$. Repeated targets simply merge intervals.
Proposition \ref{prop:directional-base-stock} then imposes $x\geq y$.

At $x=U$, the endpoint envelope with values $thU$ and
$\hat V_\downarrow^{t-1}(0,U)$ yields
\[
B_{\downarrow,t}(U,\mu)
=(c+th)U+
 \left(\frac{\hat V_\downarrow^{t-1}(0,U)}{U}-th\right)_{+}\mu.
\]
For $x\geq\chi_t^\downarrow(\mu)$, every active origin face has
$\mu\leq q_{t,j}$. Integrating the slopes
$\Phi_t-K_t\mu/q_{t,j}$ on these faces and $\Phi_t$ above the
highest target, starting at capacity, proves
\eqref{eq:downward-value} on this action interval. Evaluating it at
$x=\max\{y,\chi_t^\downarrow(\mu)\}$ proves the value formula.
Finally, if $z_n<U$, the induction hypothesis gives
$\hat V_\downarrow^n(0,U)=(c+nh)U-P(U-z_n)$, so
\[
A_t=\left(c-h-P(1-z_n/U)\right)_{+}
=K_t\left(\frac{P(z_n/U-\alpha)+b}{K_t}-1\right)_{+}.
\]
It vanishes when the highest target is below $U$; for $z_n=U$
the target remains $U$. This closes the induction.
\Halmos
\endproof

\proof{Proof of Corollary \ref{cor:supermartingale-stabilization}.}
Use $\lambda_t$ from the proof of
Proposition~\ref{prop:threshold-stabilization}, and set
$u_t=s_{t,t}^\downarrow/U$.
Before clipping at one, the highest-target recursion is
\[
f_t(u)=\frac{\Phi_{t-1}(u-\alpha)+b}{K_t},\qquad
f_t(u)-u=
\frac{(c-h)u+(b-\alpha\Phi_{t-1})(1-u)}{K_t}.
\]
Since $b-\alpha\Phi_{t-1}\geq
b-h-\alpha(c-w)>0$, $u_t$ increases from
$u_1=b/(b+h)\geq\alpha$. Here $\alpha>0$, so $K_t$ is bounded
above. The displayed increment is bounded below, before clipping,
by $\min\{c-h,b-h-\alpha(c-w)\}/\sup_{r\ge2}K_r>0$.
Consequently $u_t$ reaches one after finitely many periods.
For a target of age $\ell\geq1$, write
$r_{t,\ell}^\downarrow=s_{t,t-\ell}^\downarrow/U$ and
$r_{t,0}^\downarrow=u_t$. Its recursion is
\[
r_{t,\ell}^\downarrow
=\lambda_t(r_{t-1,\ell-1}^\downarrow-\alpha)_{+}.
\]
Since $u_t$ and $\lambda_t$ are nondecreasing, fixed-age values
are nondecreasing; positive ages persist and zeros cannot generate
older positive ages. Thus
$N_{t+1}^\downarrow\in\{N_t^\downarrow,N_t^\downarrow+1\}$.

The limiting age recursion starts at one and uses $\lambda_\infty$.
Its age-$\ell$ value is positive exactly when
\[
\lambda_\infty^{\ell-1}>
\frac{\alpha}{1-\beta\lambda_\infty}.
\]
Exactly ages $1,\ldots,N_\infty$ satisfy this test; they become
positive in finite time, while other ages remain zero by monotonicity.
The count stabilizes upon $s_{t,t-N_\infty}^\downarrow>0$, proving
$T_\downarrow^{\rm stab}$.
\Halmos
\endproof

\section{Numerical Implementation}
\label{app:history-study}
The master seed is 2026102004.

\emph{Demand and histories.}
The conditional Beta law is discretized on $\{0,\ldots,U\}$
by unbiased rounding: $z\in[k,k+1]$ rounds to $k+1$ with
probability $z-k$ and to $k$ otherwise. Integrating the Beta law
under this rule gives integer-demand mean
$(1-\delta)m+\delta p_{\mathrm B}U$.
Each of 80 replications draws operating parameters and
$p_{\mathrm B}$; these are shared across demand settings. The
central-cost sweep uses instances feasible at every tested $c$.
Independent 24-transition training episodes start at demand 10;
$N\in\{24,48,72,96\}$ uses one to four episodes without
cross-reset transitions. Training and test streams are independent.

\emph{Policies.}
M-DRO uses the analytical policy. Rolling CE optimizes assuming
constant future demand equal to the latest observation. For training
pairs $(m_i,d_i)$, let $n_{md}=\sum_i\mathbf1_{\{m_i=m,d_i=d\}}$
and $n_m=\sum_d n_{md}$. IID DP uses $p_d^{\mathrm{IID}}$ at
every state and Conditional DP uses $P_{md}^{\mathrm{cond}}$:
\[
 p_d^{\mathrm{IID}}=
 \frac{\sum_i\mathbf1_{\{d_i=d\}}+0.5/(U+1)}{N+0.5},\qquad
 P_{md}^{\mathrm{cond}}=\frac{n_{md}}{n_m}\quad(n_m>0).
\]
If $n_m=0$, set
$\mathcal S(m)=\arg\min_{r:n_r>0}|r-m|$ and
$P_{md}^{\mathrm{cond}}=\sum_{r\in\mathcal S(m)}n_{rd}/
\sum_{r\in\mathcal S(m)}n_r$, pooling ties.
Myopic Newsvendor uses Gaussian-kernel rows with bandwidth
$\ell=3$, shrinkage $s=2$, and pseudocount $0.05$:
\[
\begin{gathered}
 p_d^{\mathrm{pool}}=
 \frac{\sum_i\mathbf1_{\{d_i=d\}}+0.05}{N+0.05(U+1)},\\
 P_{md}^{\mathrm{NV}}=
 \frac{\sum_i e^{-(m-m_i)^2/(2\ell^2)}\mathbf1_{\{d_i=d\}}+
 s p_d^{\mathrm{pool}}}
 {\sum_i e^{-(m-m_i)^2/(2\ell^2)}+s}.
\end{gathered}
\]
It orders up to the smallest fitted $(b-c+\alpha w)/(b+h)$
quantile, subject to inventory. Fixed Base-stock minimizes training
cost over $\{0,0.5,\ldots,U\}$ using 24-period episodes with
initial inventory $U/2$ and demand 10, also for horizons 12 and 48.
Rolling CE and both DPs use a 0.125-unit inventory grid, linear
continuation interpolation, and exact no-order actions, without
rounding state transitions.

\emph{Evaluation.}
Each replication generates 250 test paths shared by all policies.
Common random numbers couple settings; horizons 12, 24, and 48
share path prefixes.
For Table~\ref{tab:beta-robustness}, each replication draws
$\delta\sim\mathrm{Unif}(0.075,0.30)$ and
$\kappa\sim\mathrm{Unif}(0.5,2)$ independently of costs and
$p_{\mathrm B}$. Four benchmarks fit $N=48$ historical transitions
once, without updating. High, medium, and low dispersion use
$\kappa/2$, $2\kappa$, and $4\kappa$; reversion shifts use
$\delta+0.10$ or $\delta-0.05$; joint shifts use $\delta+0.10$
with $\kappa/2$ or $2\kappa$; IID demand sets $\delta=1$.
Operating parameters, $p_{\mathrm B}$, and training histories stay
fixed within each replication; test paths use common random numbers
but may differ across shifts.
Costs and supply/recovery average 250 paths and $T$ periods
($T-1$ for nonterminal supply). Fill rate is
$1-\text{lost demand}/\text{total demand}$ (one for zero demand).
Pointwise intervals use 80 replications and Student's $t$ with
79 degrees of freedom. Savings
$100(1-\overline C^{\mathrm{M\text{-}DRO}}/\overline C^\pi)$
use the paired delta method; absolute differences use paired
replications. All quantities are per period.

%% file: MSOM_submission.bbl
\begin{thebibliography}{57}
\providecommand{\natexlab}[1]{#1}
\providecommand{\url}[1]{\texttt{#1}}
\providecommand{\urlprefix}{URL }

\bibitem[{Adida \protect\BIBand{} Perakis(2006)}]{adida2006robust}
Adida E, Perakis G (2006) A robust optimization approach to dynamic pricing and
  inventory control with no backorders. \emph{Mathematical Programming}
  107(1):97--129.

\bibitem[{Ahmed et~al.(2007)Ahmed, {\c{C}}akmak, \protect\BIBand{}
  Shapiro}]{ahmed2007coherent}
Ahmed S, {\c{C}}akmak U, Shapiro A (2007) Coherent risk measures in inventory
  problems. \emph{European Journal of Operational Research} 182(1):226--238.

\bibitem[{Asadi \protect\BIBand{} Nurre~Pinkley(2022)}]{asadi2022monotone}
Asadi A, Nurre~Pinkley S (2022) A monotone approximate dynamic programming
  approach for the stochastic scheduling, allocation, and inventory
  replenishment problem: Applications to drone and electric vehicle battery
  swap stations. \emph{Transportation Science} 56(4):1085--1110.

\bibitem[{Asadi \protect\BIBand{} Pinkley(2021)}]{asadi2021stochastic}
Asadi A, Pinkley SN (2021) A stochastic scheduling, allocation, and inventory
  replenishment problem for battery swap stations. \emph{Transportation
  Research Part E: Logistics and Transportation Review} 146:102212.

\bibitem[{Badinelli(1990)}]{badinelli1990inventory}
Badinelli RD (1990) The inventory costs of common mis-specification of
  demand-forecasting models. \emph{International Journal of Production
  Research} 28(12):2321--2340.

\bibitem[{Bertsimas et~al.(2010)Bertsimas, Doan, Natarajan, \protect\BIBand{}
  Teo}]{bertsimas2010models}
Bertsimas D, Doan XV, Natarajan K, Teo CP (2010) Models for minimax stochastic
  linear optimization problems with risk aversion. \emph{Mathematics of
  Operations Research} 35(3):580--602.

\bibitem[{Chen et~al.(2000)Chen, Drezner, Ryan, \protect\BIBand{}
  Simchi-Levi}]{chen2000quantifying}
Chen F, Drezner Z, Ryan JK, Simchi-Levi D (2000) Quantifying the bullwhip
  effect in a simple supply chain: The impact of forecasting, lead times, and
  information. \emph{Management Science} 46(3):436--443.

\bibitem[{Chen et~al.(2024)Chen, Xia, Shao, Gu, Wu, Wu, \protect\BIBand{}
  Wang}]{chen2024charging}
Chen Z, Xia S, Shao Z, Gu Z, Wu Z, Wu L, Wang X (2024) Charging or swapping? a
  study on the private consumers’ acceptance of the electricity replenishment
  mode of new energy vehicles. \emph{Transportation Research Part F: Traffic
  Psychology and Behaviour} 107:903--923.

\bibitem[{Dana~Jr \protect\BIBand{} Petruzzi(2001)}]{dana2001note}
Dana~Jr JD, Petruzzi NC (2001) Note: The newsvendor model with endogenous
  demand. \emph{Management Science} 47(11):1488--1497.

\bibitem[{Dehaybe et~al.(2024)Dehaybe, Catanzaro, \protect\BIBand{}
  Chevalier}]{dehaybe2024deep}
Dehaybe H, Catanzaro D, Chevalier P (2024) Deep reinforcement learning for
  inventory optimization with non-stationary uncertain demand. \emph{European
  Journal of Operational Research} 314(2):433--445.

\bibitem[{Delage \protect\BIBand{} Ye(2010)}]{delage2010distributionally}
Delage E, Ye Y (2010) Distributionally robust optimization under moment
  uncertainty with application to data-driven problems. \emph{Operations
  Research} 58(3):595--612.

\bibitem[{Dong \protect\BIBand{} Lee(2003)}]{dong2003optimal}
Dong L, Lee HL (2003) Optimal policies and approximations for a serial
  multiechelon inventory system with time-correlated demand. \emph{Operations
  Research} 51(6):969--980.

\bibitem[{Fildes \protect\BIBand{} Makridakis(1995)}]{fildes1995impact}
Fildes R, Makridakis S (1995) The impact of empirical accuracy studies on time
  series analysis and forecasting. \emph{International Statistical Review/Revue
  Internationale de Statistique} 289--308.

\bibitem[{Fildes et~al.(2008)Fildes, Nikolopoulos, Crone, \protect\BIBand{}
  Syntetos}]{fildes2008forecasting}
Fildes R, Nikolopoulos K, Crone SF, Syntetos AA (2008) Forecasting and
  operational research: a review. \emph{Journal of the Operational Research
  Society} 59(9):1150--1172.

\bibitem[{Gallego \protect\BIBand{} {\"O}zer(2001)}]{gallego2001integrating}
Gallego G, {\"O}zer {\"O} (2001) Integrating replenishment decisions with
  advance demand information. \emph{Management Science} 47(10):1344--1360.

\bibitem[{Gao et~al.(2024)Gao, Chen, \protect\BIBand{}
  Kleywegt}]{gao2024wasserstein}
Gao R, Chen X, Kleywegt AJ (2024) Wasserstein distributionally robust
  optimization and variation regularization. \emph{Operations Research}
  72(3):1177--1191.

\bibitem[{Hanasusanto et~al.(2015)Hanasusanto, Kuhn, Wallace, \protect\BIBand{}
  Zymler}]{hanasusanto2015distributionally}
Hanasusanto GA, Kuhn D, Wallace SW, Zymler S (2015) Distributionally robust
  multi-item newsvendor problems with multimodal demand distributions.
  \emph{Mathematical Programming} 152(1):1--32.

\bibitem[{He et~al.(2022)He, Zhu, Li, Chen, \protect\BIBand{}
  Wu}]{he2022sizing}
He C, Zhu J, Li S, Chen Z, Wu W (2022) Sizing and locating planning of ev
  centralized-battery-charging-station considering battery logistics system.
  \emph{IEEE Transactions on Industry Applications} 58(4):5184--5197.

\bibitem[{{International Energy Agency}(2025)}]{IEA_GlobalEVOutlook2025}
{International Energy Agency} (2025) \emph{Global EV Outlook 2025} (IEA,
  Paris), \urlprefix\url{https://www.iea.org/reports/global-ev-outlook-2025},
  accessed: 2025-12-2.

\bibitem[{Iyengar(2005)}]{iyengar2005robust}
Iyengar GN (2005) Robust dynamic programming. \emph{Mathematics of Operations
  Research} 30(2):257--280.

\bibitem[{Kang et~al.(2015)Kang, Brisimi, \protect\BIBand{}
  Paschalidis}]{kang2015distribution}
Kang SC, Brisimi TS, Paschalidis IC (2015) Distribution-dependent robust linear
  optimization with applications to inventory control. \emph{Annals of
  Operations Research} 231(1):229--263.

\bibitem[{Kim \protect\BIBand{} Do~Chung(2024)}]{kim2024data}
Kim YG, Do~Chung B (2024) Data-driven wasserstein distributionally robust
  dual-sourcing inventory model under uncertain demand. \emph{Omega}
  127:103112.

\bibitem[{Klabjan et~al.(2013)Klabjan, Simchi-Levi, \protect\BIBand{}
  Song}]{klabjan2013robust}
Klabjan D, Simchi-Levi D, Song M (2013) Robust stochastic lot-sizing by means
  of histograms. \emph{Production and Operations Management} 22(3):691--710.

\bibitem[{Li et~al.(2022)Li, Wang, Li, Yi, \protect\BIBand{}
  Qi}]{li2022battery}
Li C, Wang N, Li W, Yi Q, Qi D (2022) A battery centralized scheduling strategy
  for battery swapping of electric vehicles. \emph{Journal of Energy Storage}
  51:104327.

\bibitem[{Liang et~al.(2022)Liang, Ding, Zhao, \protect\BIBand{}
  Lee}]{liang2022real}
Liang Y, Ding Z, Zhao T, Lee WJ (2022) Real-time operation management for
  battery swapping-charging system via multi-agent deep reinforcement learning.
  \emph{IEEE Transactions on Smart Grid} 14(1):559--571.

\bibitem[{Liu et~al.(2018)Liu, Zhao, Yao, Soh, \protect\BIBand{}
  Wang}]{liu2018distributed}
Liu X, Zhao T, Yao S, Soh CB, Wang P (2018) Distributed operation management of
  battery swapping-charging systems. \emph{IEEE Transactions on Smart Grid}
  10(5):5320--5333.

\bibitem[{Maggioni et~al.(2025)Maggioni, Dabbene, \protect\BIBand{}
  Pflug}]{maggioni2025sampling}
Maggioni F, Dabbene F, Pflug GC (2025) Sampling methods for multi-stage robust
  optimization problems. \emph{Annals of Operations Research} 1--39.

\bibitem[{Mak et~al.(2013)Mak, Rong, \protect\BIBand{}
  Shen}]{mak2013infrastructure}
Mak HY, Rong Y, Shen ZJM (2013) Infrastructure planning for electric vehicles
  with battery swapping. \emph{Management Science} 59(7):1557--1575.

\bibitem[{Ni et~al.(2020)Ni, Sun, Tan, \protect\BIBand{}
  Tsang}]{ni2020inventory}
Ni L, Sun B, Tan X, Tsang DH (2020) Inventory planning and real-time routing
  for network of electric vehicle battery-swapping stations. \emph{IEEE
  Transactions on Transportation Electrification} 7(2):542--553.

\bibitem[{Patel et~al.(2024)Patel, Arora, Singh, Mahapatra, Chaturvedi,
  \protect\BIBand{} Saini}]{patel2024impact}
Patel M, Arora P, Singh R, Mahapatra D, Chaturvedi V, Saini SK (2024) Impact of
  battery swapping in the passenger sector: Ev adoption, emissions, and energy
  mix. \emph{Energy} 298:131393.

\bibitem[{Pichler \protect\BIBand{} Shapiro(2021)}]{pichler2021mathematical}
Pichler A, Shapiro A (2021) Mathematical foundations of distributionally robust
  multistage optimization. \emph{SIAM Journal on Optimization}
  31(4):3044--3067.

\bibitem[{Qi et~al.(2023)Qi, Zhang, \protect\BIBand{} Zhang}]{qi2023scaling}
Qi W, Zhang Y, Zhang N (2023) Scaling up electric-vehicle battery swapping
  services in cities: A joint location and repairable-inventory model.
  \emph{Management Science} 69(11):6855--6875.

\bibitem[{Qiu et~al.(2021)Qiu, Sun, \protect\BIBand{}
  Sun}]{qiu2021distributionally}
Qiu R, Sun Y, Sun M (2021) A distributionally robust optimization approach for
  multi-product inventory decisions with budget constraint and demand and yield
  uncertainties. \emph{Computers \& Operations Research} 126:105081.

\bibitem[{Raad et~al.(2025)Raad, Rajendran, Srinivas, Khowal, Abam,
  \protect\BIBand{} Lewis}]{raad2025hybrid}
Raad NG, Rajendran S, Srinivas S, Khowal A, Abam B, Lewis E (2025) Hybrid and
  distributionally robust optimization approaches for hospital blood bank
  inventory management under demand uncertainty. \emph{Journal of the
  Operational Research Society} 1--24.

\bibitem[{Sayman \protect\BIBand{}
  {\"O}nc{\"u}ler(2009)}]{sayman2009investigation}
Sayman S, {\"O}nc{\"u}ler A (2009) An investigation of time inconsistency.
  \emph{Management Science} 55(3):470--482.

\bibitem[{Scarf(1959)}]{scarf1959bayes}
Scarf H (1959) Bayes solutions of the statistical inventory problem. \emph{The
  Annals of Mathematical Statistics} 30(2):490--508.

\bibitem[{Scarf(1960)}]{scarf1960some}
Scarf HE (1960) Some remarks on bayes solutions to the inventory problem.
  \emph{Naval Research Logistics Quarterly} 7(4):591--596.

\bibitem[{Scarf et~al.(1957)Scarf, Arrow, \protect\BIBand{}
  Karlin}]{scarf1957min}
Scarf HE, Arrow K, Karlin S (1957) A min-max solution of an inventory problem.
  Technical report, Rand Corporation Santa Monica.

\bibitem[{Schneider et~al.(2018)Schneider, Thonemann, \protect\BIBand{}
  Klabjan}]{schneider2018optimization}
Schneider F, Thonemann UW, Klabjan D (2018) Optimization of battery charging
  and purchasing at electric vehicle battery swap stations.
  \emph{Transportation Science} 52(5):1211--1234.

\bibitem[{Shaker et~al.(2023)Shaker, Farzin, \protect\BIBand{}
  Mashhour}]{shaker2023joint}
Shaker MH, Farzin H, Mashhour E (2023) Joint planning of electric vehicle
  battery swapping stations and distribution grid with centralized charging.
  \emph{Journal of Energy Storage} 58:106455.

\bibitem[{Shanghai(2022)}]{Shanghai2022}
Shanghai (2022) Implementation opinions on further promoting the construction
  of charging and swapping infrastructure in this city.
  \urlprefix\url{https://www.shanghai.gov.cn/xnysj/20250610/ed23a7c382c1498abc0656b768cefd9e.html},
  accessed: 2025-12-2.

\bibitem[{Shapiro(2012)}]{shapiro2012minimax}
Shapiro A (2012) Minimax and risk averse multistage stochastic programming.
  \emph{European Journal of Operational Research} 219(3):719--726.

\bibitem[{Shapiro(2016)}]{shapiro2016rectangular}
Shapiro A (2016) Rectangular sets of probability measures. \emph{Operations
  Research} 64(2):528--541.

\bibitem[{Shapiro \protect\BIBand{} Homem-de Mello(2000)}]{shapiro2000rate}
Shapiro A, Homem-de Mello T (2000) On the rate of convergence of optimal
  solutions of monte carlo approximations of stochastic programs. \emph{SIAM
  Journal on Optimization} 11(1):70--86.

\bibitem[{Shapiro \protect\BIBand{} Kleywegt(2002)}]{shapiro2002minimax}
Shapiro A, Kleywegt A (2002) Minimax analysis of stochastic problems.
  \emph{Optimization Methods and Software} 17(3):523--542.

\bibitem[{Shenzhen(2023)}]{Shenzhen2023}
Shenzhen (2023) Shenzhen municipality's management measures for charging and
  swapping facilities for new energy vehicles.
  \urlprefix\url{https://www.sz.gov.cn/zfgb/2023/gb1307/content/post_10906702.html},
  accessed: 2025-12-2.

\bibitem[{Sun et~al.(2019{\natexlab{a}})Sun, Sun, Tsang, \protect\BIBand{}
  Whitt}]{sun2019optimal}
Sun B, Sun X, Tsang DH, Whitt W (2019{\natexlab{a}}) Optimal battery purchasing
  and charging strategy at electric vehicle battery swap stations.
  \emph{European Journal of Operational Research} 279(2):524--539.

\bibitem[{Sun et~al.(2019{\natexlab{b}})Sun, Yang, \protect\BIBand{}
  Yang}]{sun2019robust}
Sun H, Yang J, Yang C (2019{\natexlab{b}}) A robust optimization approach to
  multi-interval location-inventory and recharging planning for electric
  vehicles. \emph{Omega} 86:59--75.

\bibitem[{Wang et~al.(2012)Wang, Atasu, \protect\BIBand{}
  Kurtulu{\c{s}}}]{wang2012multiordering}
Wang T, Atasu A, Kurtulu{\c{s}} M (2012) A multiordering newsvendor model with
  dynamic forecast evolution. \emph{Manufacturing \& Service Operations
  Management} 14(3):472--484.

\bibitem[{Wu(2021)}]{wu2021survey}
Wu H (2021) A survey of battery swapping stations for electric vehicles:
  Operation modes and decision scenarios. \emph{IEEE Transactions on
  Intelligent Transportation Systems} 23(8):10163--10185.

\bibitem[{Xie et~al.(2024)Xie, Dai, \protect\BIBand{} Pei}]{xie2024empowering}
Xie X, Dai X, Pei Z (2024) Empowering the capillary of the urban daily commute:
  Battery deployment analysis for the locker-based e-bike battery swapping.
  \emph{Transportation Science} 58(1):176--197.

\bibitem[{Xin \protect\BIBand{} Goldberg(2021)}]{xin2021time}
Xin L, Goldberg DA (2021) Time (in) consistency of multistage distributionally
  robust inventory models with moment constraints. \emph{European Journal of
  Operational Research} 289(3):1127--1141.

\bibitem[{Xin \protect\BIBand{} Goldberg(2022)}]{xin2022distributionally}
Xin L, Goldberg DA (2022) Distributionally robust inventory control when demand
  is a martingale. \emph{Mathematics of Operations Research} 47(3):2387--2414.

\bibitem[{Yin et~al.(2024)Yin, Wang, Chu, \protect\BIBand{}
  Wang}]{yin2024distributionally}
Yin Y, Wang J, Chu F, Wang D (2024) Distributionally robust multi-period
  humanitarian relief network design integrating facility location, supply
  inventory and allocation, and evacuation planning. \emph{International
  Journal of Production Research} 62(1-2):45--70.

\bibitem[{Zhang et~al.(2025{\natexlab{a}})Zhang, Turan, Sarker,
  \protect\BIBand{} Essam}]{zhang2025robusta}
Zhang D, Turan HH, Sarker R, Essam D (2025{\natexlab{a}}) Robust optimization
  approaches in inventory management: Part a—the survey. \emph{IISE
  Transactions} 57(7):818--844.

\bibitem[{Zhang et~al.(2025{\natexlab{b}})Zhang, Turan, Sarker,
  \protect\BIBand{} Essam}]{zhang2025robustb}
Zhang D, Turan HH, Sarker R, Essam D (2025{\natexlab{b}}) Robust optimization
  approaches in inventory management: Part b—the comparative study.
  \emph{IISE Transactions} 57(7):845--871.

\bibitem[{Zhang et~al.(2025{\natexlab{c}})Zhang, Lv, Kuai, \protect\BIBand{}
  Feng}]{zhang2025battery}
Zhang F, Lv H, Kuai C, Feng T (2025{\natexlab{c}}) The battery-swapping
  revolution: Exploring user preferences in electric micro-mobility sector.
  \emph{Transportation Research Part A: Policy and Practice} 194:104416.

\end{thebibliography}
